\documentstyle[amssymb,12pt]{amsart}

\newif\ifTwelve

\Twelvetrue

\DeclareFontShape{OT1}{cmr}{m}{sc}
 {<5><6><7><8>cmcsc8 <9>cmcsc9%
 <10><10.95><12><14.4><17.28><20.74><24.88> cmcsc10}{}

\ifTwelve  \else  \fi

\ifx\mathrm\undefined\let\mathrm\bf\fi
\ifx\mathbf\undefined\let\mathbf\bf\fi

\let\leq\leqslant \let\geq\geqslant

\let\Box\square
\let\al\alpha
\let\bt\beta
\let\gm\gamma 
\let\dl\delta \let\Dl\Delta
\let\epe\epsilon \let\eps\varepsilon \let\epsilon\eps

\let\la\lambda 
\let\om\omega 
 \let\phi\varphi

\newcommand{\Z}{{\Bbb Z}}

\newcommand{\C}{{\Bbb C}}

\newcommand{\Ref}[1]{{$($\ref{#1}$)$}}
\newcommand{\bean}{\begin{eqnarray}}
\newcommand{\eean}{\end{eqnarray}}
\newcommand{\be}{\begin{displaymath}}
\newcommand{\ee}{\end{displaymath}}
\newcommand{\bea}{\begin{eqnarray*}}
\newcommand{\eea}{\end{eqnarray*}}
\newcommand{\g}{{{\frak g}\,}}

\newcommand{\n}{{{\frak n}}}
\newcommand{\h}{{{\frak h\,}}}
\newcommand{\Id}{{\operatorname{Id}}}

\newcommand{\Ker}{{\operatorname{Ker}}}

\newcommand{\T}{\!\otimes\!}

\newcommand{\vs}{\vspace{.5\baselineskip}}

\newenvironment{proof}{\noindent{\it Proof\/}:}{$\;\Box$\par\vs}

\newtheorem%
{thm}{Theorem}%[section]
\newtheorem%
{proposition}[thm]{Proposition}
\newtheorem%
{lemma}[thm]{Lemma}
\newtheorem%
{lemmadef}[thm]{Lemma-Definition}
\newtheorem%
{corollary}[thm]{Corollary}
\newtheorem%
{conjecture}[thm]{Conjecture}

\newcommand{\End}{{\operatorname{End\,}}}
\newcommand{\Hom}{{\operatorname{Hom\,}}}

\newcommand{\Rm}{{{\cal R}}}
\newcommand{\V}{{\cal {V} }}
\newcommand{\Oc}{{\cal O}}

\newcommand{\wo}{{\widetilde \otimes}}
\newcommand{\bo}{{\bar {\otimes}}}

\newcommand{\E}{{\bar {E}\,}}
\newcommand{\Ti}{{\cal {T}}}
\newcommand{\Mm}{{M_{T}\,}}
\newcommand{\ta}{{\bar {\tau}\,}}
\newcommand{\Tb}{{\tilde { \otimes}\,}}
\newcommand{\Lb}{{\bar {L}\,}}

\newcommand{\tK}{{\tilde {K}\,}}
\newcommand{\TB}{{\bar { \otimes}\,}}
\newcommand{\lt}{{\tilde { \la }}}
\newcommand{\Repq}{{Rep_f(F_q)}}

\def\nli#1{\vsk#1>\noindent\ignorespaces}
\def\nl{\nli0}
\def\vsk#1>{\vskip#1\baselineskip}
\def\ftext#1{{\let\thefootnote\relax\footnotetext{\vsk-.8>\noindent #1}}}

\title[Exchange Dynamical Quantum Groups]
{Exchange Dynamical Quantum Groups }

\author[P.Etingof and A. Varchenko]
{P.\ Etingof $^{\,\star}$  and
A.\ Varchenko$^{\,\diamond}$}

\begin{document}
\maketitle

\begin{center}
{\it
$^\star$Department of Mathematics, Harvard University, Cambridge, MA 02138,
$^\diamond$Department of Mathematics, University of North Carolina
at Chapel Hill

Chapel Hill, NC 27599-3250, USA}
\end{center}

\ftext{\small{\normalsize\sl $^\star$Email\/{\rm etingof@@math.harvard.edu}}\nl
\hphantom{$*$}Supported in part by NSF grant  DMS 9700477 \nli{.2}
%{\normalsize\sl $^*$Email\/{\rm:} vt@@math.sci.osaka-u.ac.jp}\nl
%\hphantom{$*$}On leave of absense from St.\,Petersburg Branch of
%Steklov Mathematical Institute\nli{.2}
{\normalsize\sl $^\diamond$Email\/{\rm av@@math.unc.edu}}\nl
\hphantom{$*$} Supported in part by NSF grant DMS-9501290}

%{\rm mukhin@@math.unc.edu,
%av@@math.unc.edu}}

\vsk1.5>
\centerline{January, 1998}
\vsk1.8>

\begin{abstract}

%Dynamical quantum  groups were introduced by Felder. In our previous paper 
%q-alg 9708015 we recalled the notion of a Hopf algebroid and showed that 
%dynamical quantum groups can be viewed as Hopf algebroids. We also inroduced 
%and studied Hopf algebroids associated to quantum dynamical R-matrices of 
%Hecke type. In this paper, f
For any simple Lie algebra $\g$ and any complex 
number $q$ which is not zero or a nontrivial root of unity,
%but may be equal to 1 
we construct a dynamical quantum group (Hopf algebroid), whose 
representation theory is essentially the same as the representation theory 
of the quantum group $U_q(\g)$. This dynamical quantum group
 is obtained from the 
fusion and exchange relations between intertwining operators in representation 
theory of $U_q(\g)$, and is an
% seems to be an adequate 
algebraic structure 
standing behind these relations. 

\end{abstract}
\vsk>
\vsk0>

\thispagestyle{empty}

\section{Introduction}

One of the most important equations in statistical mechanics is the so-called 
Star-Trangle relation, introduced by Baxter. In 1994, G.Felder \cite{F} 
suggested to write this relation in the form of the quantum dynamical 
Yang-Baxter equation (QDYB)
(which previously appeared in some form in physical 
literature), and proposed the concept of a quantum group associated to a 
solution of this equation. He also considered the quasiclassical
limit of this equation, and showed that a solution of the classical
dynamical Yang-Baxter 
equation (CDYB)
appears naturally on the right hand side of the 
Knizhnik-Zamolodchikov-Bernard equations for conformal blocks on an 
elliptic curve. Since then, this theory has found many applications
in the theory of integrable systems. 

In [EV1], we proposed a geometric interpretation of the 
CDYB equation without spectral parameter. Namely, we  assigned to any solution 
of this equation, whose symmetric part is invariant, a certain Poisson 
groupoid. This construction is a generalization of Drinfeld's construction 
which assigns a Poisson-Lie group
 to any solution of the usual classical Yang-Baxter equation,
with symmetric invariant part. 
We also classified such solutions for simple Lie algebras 
and showed that there are two classes of solutions 
(without spectral parameter)-- rational 
and trigonometric. 

In [EV2], we quantized the results of [EV1] and presented 
a ``noncommutative geometric'' interpretation of solutions
of the QDYB equation without spectral parameter. Namely, we assigned to any solution 
of the QDYB equation, satisfying a special ``Hecke type'' condition,
a certain dynamical quantum group (Hopf algebroid). 
This construction is a generalization of 
the Faddeev-Reshetikhin-Takhtajan-Sklyanin construction 
which assigns to any solution of the usual quantum Yang-Baxter equation 
of Hecke type a quantum group (Hopf algebra), 
defined by the so called $RTT=TTR$ relations. 
We also classified the Hecke type solutions  of the QDYB equation
and showed that, 
like in the classical case,  there are two classes of solutions 
(without spectral parameter)-- rational 
and trigonometric. 

The solutions of the QDYB equation
and corresponding dynamical quantum groups from [EV2] provide 
quantizations of the solutions of the CDYB equation and Poisson groupoids from [EV1], 
but only 
for the Lie algebra $\g=sl_N$. For other simple Lie algebras, especially 
for exceptional ones, one needs to use a different method to 
quantize the Poisson groupoids from [EV1]; this method has to be applicable 
to any simple Lie algebra $\g$ and should not use its particular 
matrix realization. 

Such a method is suggested in the present paper. Namely, it turns out 
that with any simple complex Lie group $G$ (with the Lie algebra $\g$)
and a nonzero complex number $q$ 
(which is not a nontrivial root of unity but may be equal to 1) one can 
associate a Hopf algebroid $F_q(G)$, which is a quantization 
of the Poisson groupoid associated with the simple Lie algebra $\g$ in [EV1].
More precisely, 
the case $q=1$ corresponds to the Poisson groupoid 
defined by the rational solution, and $q\ne 1$ 
corresponds to the Poisson groupoid 
defined by the trigonometric solution. 

The Hopf algebroid $F_q(G)$ is constructed from representation theory of 
$G$. Namely, the structure constants of the multiplication in $F_q(G)$ are 
obtained from the structure constants of the multiplication (fusion)  
of intertwining operators between a Verma module over $U_q(\g)$
and a tensor product of a Verma module with a finite dimensional module over  $U_q(\g)$.
These structure constants 
have been known for a long time under various names
(Racah coefficients, Wigner 6j symbols) and play an important role 
in quantum physics. 

The commutation relations between generators of $F_q(G)$
are defined by certain dynamical R-matrices, which satisfy the 
quantum dynamical Yang-Baxter equations. These R-matrices are exactly 
the matrices which arise in commutation (=exchange) relations between 
intertwining operators and are therefore called the exchange matrices. 
This makes it natural to call the Hopf algebroids $F_q(G)$
{\bf the exchange dynamical quantum groups.}

The results of this paper demonstrate how to use representation 
theory to construct new quantum groups, and conversely, how
the multiplication of intertwining operators,
one of the main structures in representation theory,  
is controlled by a dynamical quantum group. 

We note that the main idea of this paper (to use 
commutation relations between intertwining operators
to obtain new quantum groups) was inspired by the pioneering paper [FR]. 

Let us briefly describe the contents of the paper. 

In Chapter 2 we introduce, for any polarized Hopf algebra, 
the fusion and exchange matrices $J(\la),R(\la)$ and 
consider their main properties. 

In Chapter 3 we recall the notion of an $H$-Hopf algebroid
and its dynamical representations introduced in [EV2]
(where $H$ is a commutative, cocommutative Hopf algebra). 

In Chapter 4 we construct Hopf algebriods defined by 
the fusion and exchange matrices. 

In Chapter 5 we specialize our construction to the case 
of simple Lie groups and quantum groups, and construct 
the Hopf algebroids $F_q(G)$. 

In Chapter 6 we compute the exchange R-matrix for the vector representation of 
$G=GL(N)$, and show that $F_q(GL(N))$ is isomorphic 
to the Hopf algebroid $A_R$ defined by the trigonometric solution $R$
of the quantum dynamical Yang-Baxter equation in [EV2]. 

In Chapter 7 we consider the representation theory of $F_q(G)$
for a simple complex group $G$ 
and show that its category of rational finite dimensional dynamical 
representations contains the category of finite dimensional representations 
of $U_q(\g)$ as a full subcategory. 
In the next paper we plan to show that 
these categories are actually the same. 

In Chapter 8, we describe the precise connection between 
fusion and exchange matrices (for $sl(2)$) and the classical and quantum 
6j-symbols. 

In Chapter 9, we show that the universal fusion matrix
$J(\l)$ satisfies the defining property of the quasi-Hopf twist  
discovered in [A]. In particular, 
this shows that our fusion matrices are the same as the 
quasi-Hopf twists introduced in [A]. 

In a subsequent paper, we plan to consider the analogue of this theory 
for affine and quantum affine algebras. This will help one to understand
better the monodromy of classical and quantum Knizhnik-Zamolodchikov equations 
following the ideas of [FR] and [TV1-2],[FTV].

In conclusion we would like to mention the paper [BBB], in which 
another algebraic interprtetation of the QDYB
equation was given (via quasi-Hopf algebras), and a version 
of our main construction for the Lie algebra $sl_2$ was presented. 
See also [JKOS] where  the approach of [BBB] was generalized to an arbitrary 
Kac-Moody algebra. We would also like to point out the recent paper 
[Xu], where the relationship between the quasi-Hopf algebra and Hopf 
algebroid interpretation of the quantum dynamical Yang-Baxter 
equation is explained. 

{\bf Acknowledgments.} We are grateful to the referee for several 
interesting questions, the answers to which have enriched the paper. 

\section{ Exchange Matrices }\label{def}

\subsection{Polarized Hopf Algebras  }\label{pha}

{\it A polarized Hopf algebra} is a Hopf algebra $A$ over $\C$
 with the following properties:
\begin{enumerate}
\item[I.] The algebra is $\Z$-graded, $A=\oplus_{k=-\infty}^\infty
A[k]$.
\item[II.] The algebra is polarized. Namely, there exist graded subalgebras
$A_0, A_+, A_-$ such that the multiplication maps $A_+\T A_0\T A_- \to A$
and $A_-\T A_0\T A_+ \to A$ are isomorphisms of vector spaces.
We also assume that the graded  components of $A_+$ and $A_-$ (not of $A_0$)
are finite dimensional.
\item[III.] Let $\epe : A \to \C $ be the counit, then 
Ker $\epe \,\cap A_+$ has only elements of positive degree,
Ker $\epe \,\cap A_-$ has only elements of negative degree,
$A_0$ has only elements of zero degree.

\item[IV.] The algebra $A_0$ is a commutative cocommutative finitely generated
Hopf algebra.
\item [V.] $ A_0 A_+$  and 
$A_0 A_-$ are Hopf subalgebras of $A$.
\end{enumerate}

Let $T=$Spec $A_0$. Since $A_0$ is finitely generated, 
commutative and cocommutative, 
$T$ is  a commutative affine algebraic group [M]. 

The main examples of polarized Hopf algebras are the universal enveloping algebra
$U(\g)$ of a simple Lie algebra $\g$ and the corresponding quantum group $U_q(\g)$.

If $\h \subset \g$ is a Cartan subalgebra, $\g= \g_- \oplus \h \oplus \g_+$ 
a  polarization and $A=U(\g)$, then $A_-=U(\g_-)$, $A_0=U(\h)$, 
$A_+=U(\g_+)$, and $T=\C^n$, where $n=$dim $\h$. 
If $A=U_q(\g)$, then $A_-=U_q(\g_-)$, $A_0=U_q(\h)$, 
$A_+=U_q(\g_+)$, and $T=(\C^*)^n$.

\subsection{ Verma modules }\label{vm}

Let $A^0_+=\Ker \,\epe\,\cap A_+$. Then $A_+=\C\cdot 1 \oplus A_+^0$.

\begin{lemma}\label{} $A_0A^0_+$ is the set of elements of $A_0A_+$ of positive degree.
Moreover, $A_0A^0_+$ is a two-sided ideal in 
$A_0A_+$.

\end{lemma}
\begin{proof} 
$A_0A^0_+$ obviously lies among the elements of positive degree. We prove the converse
statement. Let $\al = \sum a^i_0a^i_+$ be a homogeneous element of degree $N > 0$ in
$A_0A_+$, where $a^i_0\in A_0$ and $a_+^i \in A_+$ has degree $N$. Since
$A_+=\C\cdot 1 \oplus A_+^0$, we have $a_+^i \in A^0_+$. The second statement of the Lemma
obviously follows from the first.

\end{proof}

\begin{corollary}\label{}
$A_0A_+/A_0A_+^0 \cong A_0$.

\end{corollary}

In fact, $A_0\to A_0A_+ \to A_0A_+/A_0A_+^0 $ is an isomorphism.
%the set of elements of degree zero, and this is $A_0$.

Let $\phi_+ : A_0A_+ \to A_0$ be the induced homomorphism, $\phi_+$ is defined by
$\phi_+(a_0)=a_0, \, \phi_+(a_+)=\epe(a_+)$ where $a_0\in A_0, a_+\in A_+$.

Let $\lambda : A_0 \to \C$ be a homomorphism, hence $\la\in T$.
Let $\chi_\la$ be this one dimensional $A_0$-module. 
Define a homomorphism $\chi_\la^+ : A_0A_+ \to \C$ by
$\chi_\la^+ =\la \phi_+$. We also denote $\chi_\la^+$ the corresponding one-dimensional
$A_0A_+$-module.
Define {\it a Verma module} $M^+_\la$  over $A$ by $M^+_\la=A\T_{A_0A_+} \chi^+_\la$.

Analogously, consider the homomorphism
$\phi_- : A_0A_- \to A_0$ and the corresponding $A_0A_-$-module
$\chi^-_\la$. Define {\it a Verma module} $M^-_\la$ over $A$ by $M^-_\la=A\T_{A_0A_-}\chi^-_\la$.

Let $\chi_\la^+=\C v_\la^+$.

\begin{lemma}\label{free}
$M^+_\la$ is a free $A_-$-module generated by $v_\la^+$.

\end{lemma}
The Lemma follows from Property II.

Similarly, $M^-_\la$ is a free $A_+$-module generated by $v_\la^-$.

Using Lemma \ref{free} induce a grading on $M^{\pm}_\la$ from $A_{\pm}$
such that the degree of $v^{\pm}_\la$ is equal to zero.

\subsection{ The Shapovalov form }\label{sf}

A polarized Hopf algebra $A$ is called {\it nondegenerate} if Verma modules $M^+_\la$
and $M^-_\la$ are irreducible for generic $\la \in T$. (This means that the modules are
irreducible for all $\la$ except a countable union of algebraic sets of lower dimension.)

Let $A$ be polarized and nondegenerate. Consider the vector space
$(M^+_\la)^*$, the restricted dual of $M^+_\la$ with respect to the grading of $A_-$.
Define an $A$-module structure on $(M^+_\la)^*$ by $\pi_{(M^+_\la)^*}(a)=
\pi_{M^+_\la}(S(a))^*$ where $S$ is the antipode in $A$.

Define $(v_\la^+)^*\in (M^+_\la)^* $ to be the only degree zero element such that
$<(v_\la^+)^*, v_\la^+>=1$.

\begin{lemma}
If $a_-\in A^0_-$, then $a_-(v_\la^+)^*=0$.
If $a_0\in A_0$, then $a_0(v_\la^+)^*=(-\la)(a_0)(v_\la^+)^*$
where $-\la$ means the inverse element in the abelian group $T$.
\end{lemma}
In fact, $a_-(v_\la^+)^*$ has degree which does not occur in $(M^+_\la)^*$.
The second statement is obvious. 

\begin{corollary} There exists a unique homomorphism of $A$-modules
$\psi_-: M_{-\la}^- \to (M_\la^+)^*$ such that $v_{-\la}^- \mapsto (v_\la^+)^*$.
\end{corollary}

Analogously, one can define a homomorphism
$\psi_+: M_{-\la}^+ \to (M_\la^-)^*$ such that $v_{-\la}^+ \mapsto (v_\la^-)^*$.

\begin{lemma}
If $A$ is nondegenerate, then $\psi_+, \psi_-$ are isomorphisms for generic $\la$.
\end{lemma}
\begin{proof} For generic $\la$, the homomorphisms $\psi_+$ and $ \psi_-$ are injective
since $M^+_{-\la}$ and $M^-_{-\la}$ are irreducible and $\psi_+, \psi_-$  differ from zero.
This implies that dim $A_+[n]\leq$ dim $ A_-[-n]$ and
dim $A_+[n]\geq$ dim $ A_-[-n]$. 
\end{proof}
\begin{corollary} 
For all $n$, dim $A_+[n]=$ dim $ A_-[-n]$.
\end{corollary}

The homomorphisms $\psi_+, \psi_-$ are called {\it
the Shapovalov forms}. They can be considered
as bilinear forms $\psi_+(\la): A_+\T A_- \to\C ,\,$
$\psi_-(\la): A_-\T A_+ \to\C $ depending on $\la$. The bilinear forms can be  defined
also by $\psi_+(\la) (a_+,a_-)= <(v_\la^+)^*, S(a_+)a_-v_\la^+>$ and
$\psi_-(\la) (a_-,a_+)= <(v_\la^-)^*, S(a_-)a_+v_\la^->$.

Choose  bases in $A_{\pm}[n]$ and compute the determinants of the Shapovalov forms,
$$
D_n^+(\la)= \text{det}\, \psi_+(\la)[n], \qquad
D_n^-(\la)= \text{det} \,\psi_-(\la)[n].
$$
The determinants of the Shapovalov forms are regular nonzero functions of $\la\in T$
defined up to multiplication by a nonzero number.

\subsection{ Intertwining operators }\label{io}
Let $A$ be polarized and nondegenerate.
Let $V$ be a $\Z$-graded $A$-module such that $V$ is a diagonalizable
$A_0$-module, $V=\oplus_{\la\in T}V[\la]$ where $a_0v=\la(a_0)v$ for all
$v\in V[\la], a_0\in A_0$.
% Assume that $V$ is a $\Z$-graded $A$-module.

\begin{thm} \label{} ${}$

\begin{enumerate}
\item[I.] Assume that $M^-_{-\mu}$ 
%(look at ${}^-$ and $-\mu$!) 
is irreducible and $V$ is bounded from above, i.e.
the graded component of $V$ corresponding a  number $N$ is equal to zero if
$N>>0$. Then $\Hom_A (M^+_\la, M^+_\mu \T V)\cong V[\la-\mu]$,
where $\la-\mu$ means the difference in the abelian group $T$.
\item[II.] Assume that $M^+_{-\mu}$ 
%(look at ${}^+$ and $-\mu$!) 
is irreducible and $V$ is bounded from below.
Then $\Hom_A (M^-_\la, M^-_\mu \T V)\cong V[\la-\mu]$.
\end{enumerate}
The isomorphism is given by
$$
\Phi \mapsto <\Phi>:=<(v^\pm_\mu)^*\T\,1,\,\Phi v_\la^\pm>.
$$
\end{thm}
\begin{proof}
First we prove a Lemma. 
Let $B$ be a $\Z$-graded Hopf algebra.
Let $U, W\,$ be $\Z$-graded $B$-modules bounded from above 
and such that all homogeneous components of $U$ are finite dimensional.
Define the space $\Hom_B (U^*, W)$ as $\Hom_B (U^*, W)=$ \newline
$\oplus_n \Hom_B (U^*, W)[n]$. 
Let $(U\T W)^{B}$ denote the subspace of invariants with respect to
$B$, i.e. the subspace of all elements $w$ such that $bw=\epe(b)w$ for any $b\in B$. 
The space $(U\T W)^B$ is $\Z$-graded, $(U\T W)^B=\oplus_n
(U\T W)^B[n]$.

\begin{lemma}
Let $w\in U\T W$. Let $\tilde w:
U^* \to W$ be defined as the composition of the two maps: 
$1\,\T\, w: U^* \to U^*\T U \T W$
and $<\cdot,\cdot >\T \,1 : U^*\T U\T W \to W$. Then $w\in (U\T W)^B[n]$ 
for some $n$ if and only if 
$\tilde w \in \Hom_B (U^*, W)[n]$. Thus the assignment $w\to \bar {w}$ is an
isomorphism of $(U\T W)^B$ and $\Hom_B(U^*,W)$.

\end{lemma}
\begin{proof}
The counit of $B$ defines a trivial
one-dimensional $B$-module $\C_B$. If $w\in U\T W$ is $B$-invariant,
then $w$ defines a homomorphism $\C_B \to U\T W$ of $B$-modules, 
and which  defines
a homomorphism $1\,\T\, w \,: \,
U^*\to U^*\T U\T W$, $u^* \mapsto u^*\T w$. Since $<\,,\,>: U^*\T U \to \C$
is a homomorphism of $B$-modules, so is the composition $\tilde w=(<\,,\,>\T\,1)\,(1\,\T
\,w).$ This proves one of the two claims of the Lemma.

Let $U\,\tilde {\T}\, U^*=\oplus_n U \,\tilde {\T}\, U^*[n]$ where
$U\,\tilde {\T}\, U^*[n]$ is the space of all elements $x$ of the form $ x=\sum_{i=1}^\infty a_i\T b_i$
such that
for any $i$ the element $a_i\T b_i \in U \T U^*$ has degree $n$, the elements 
$a_i, b_i$ are homogeneous, and
deg $a_i \to -\infty$ as $i\to \infty$.
There is a canonical $B$-homomorphism $\C_B \to U\,\tilde {\T} \,U^*$, 
$1 \mapsto \sum_i a_i \T a^*_i$ where $\{a_i\}$ is a graded basis
of $U$ and $\{a_i^*\}$ is  the dual basis of $U^*$.

Let $\tilde w : U^*\to W$ be a homogeneous $B$-homomorphism. Then
$1\,\T\,\tilde w : U \,\tilde {\T}\, U^* \to U\T W$ is a well defined $B$-homomorphism.
The composition $\C_B \to U \,\tilde {\T}\, U^* \to U\T W$ gives a $B$ invariant element $w$.

\end{proof}

Now we prove the Theorem.
Introduce $A_{\geq 0}=A_0A_+$ and $A_{\leq 0}=A_0A_-$. We have 
$\Hom _A(M^+_\la,M^+_\mu\T V)=\Hom_{A_{\geq 0}}(\chi_\la^+, M^+_\mu\T V)$.
The space $\Hom_{A_{\geq 0}}(\chi_\la^+, M^+_\mu\T V)$ 
can be described as the space $X$
of all $w\in M^+_\mu\T V$ such that the
$A_{\geq 0}$-submodule of $M^+_\mu\T V$ generated by $w$ is isomorphic to $\chi^+_\la$.
After tensoring with $\chi^+_{-\la}$ this submodule gives a trivial module.
Thus the space
$X$ is isomorphic to the space $ (M^+_\mu\T V\T \chi^+_{-\la})^{A_{\geq 0}}$.
According to the Lemma, the space $X$ is isomorphic to 
$\Hom_{A_{\geq 0}}((M^+_\mu)^*, V\T \chi^+_{-\la})$.
This space is isomorphic to
$\Hom_{A_{\geq 0}}(M^-_{-\mu}, V\T \chi^+_{-\la})$
since $M^-_{-\mu}$ is irreducible. Now
$\Hom_{A_{\geq 0}}(M^-_{-\mu}, V\T \chi^+_{-\la})\cong 
\Hom_{A_{0}}(\chi_{-\mu}, V\T \chi_{-\la})\cong
\Hom_{A_{0}}(\chi_{\la-\mu}, V)\cong V[\la-\mu]$. The Theorem is proved.

\end{proof}

Let $V$ be bounded from above. Let $v\in V[\la-\mu]$.
Denote $\Phi_\la^v:M^+_\la \to M^+_\mu\T V$ the intertwining operator
such that $<\Phi^v_\la>=v$.

Define $\Phi^v(\la): A_-\to A_-\T V$ as
the operator obtained from $\Phi^v_\la$ after identification of $A_-$
with $M^+_\la$ and $ M^+_\mu$. Then $\Phi^v(\la)$ is a rational function,
i.e. for any homogeneous $a_-\in A_-, a_-^*\in A_-^*$, and $f\in V^*$,
the scalar function $(f^*\T a^*_-)\Phi^v(\la)a_-$ is a rational function.

\subsection{ A quasitriangular structure and dynamical R-matrices }\label{qs}
Let $A_{\geq n}=\oplus_{j\geq n} A[j]$. Introduce a system of left
$A$-ideals $I_n= A\cdot A_{\geq n}$. 

Introduce a tensor product  $A \hat {\T} A$ $ =$ $ \oplus _{i\in \Z} ( A \, \hat {\T} \,A )[i]$
as follows. Let
$(A\, \hat {\T} \, A )[i]$ be the projective limit of $(A/I_n\T A/I_n)[i]$ as 
$n\to\infty$, that is $(A\, \hat {\T} \, A )[i]$
consist of elements of the form
$a=\sum_{k=1}^\infty a_k\T a_k'$ such that
\begin{enumerate}
\item[I.] For each $k$ there is $j$ such that 
$a_k\T a_k'\in A[j]\T A[i-j]$.
\item[II.] For each $n$ there is only finitely many $k$ such that 
$a_k\T a_k'$ does not belong to $A\T I_n + I_n\T A$.
\end{enumerate}

\begin{lemma}
$A\, \hat {\T}\, A$ is an algebra. $\square$
\end{lemma}

Similarly we can define $A^{\hat {\T}n}$ for any $n$.

Consider the category $\cal O$ of graded $A$-modules bounded from above and
diagonalizable over $A_0$. Let $V, W \in \cal O$, then  $A \hat {\T} A$ acts on
$V\T W$. Similarly, for any $n$ one can define an action of the 
algebra $ A^{\hat {\T} n}$
in a tensor product of $n\,$ $A$-modules from the category $\cal O$.

%Let $\Rm\in A\,\hat {\T} \,A$ and $V, W $ bounded from above, then $\Rm$ defines an endomorphism of
%a vector space $V\T W$.

An element $\Rm\in A\,\hat {\T}\, A$ is called {\it a quasitriangular structure   } (QTS) if
\begin{enumerate}
\item[I.] $\Rm$ is invertible in $A\,\hat{\T}\, A$.
\item[II.] $\Rm\Dl(a)=\Dl^{op}(a)\Rm$. 
\item[III.] $(\Dl\T \,1)\Rm=\Rm^{13}\Rm^{23}$ and
$(1\,\T \Dl)\Rm=\Rm^{13}\Rm^{12}$. 
\end{enumerate}

Consider the category $\cal O$ of graded $A$-modules bounded from above and
diagonalizable over $A_0$. The category $\cal O$ is a braided tensor category
with the braiding  equal to $P\Rm$.

Let $A$ be a nondegenerate polarized Hopf algebra with a QTS $\Rm$.
Let $V, W \in \cal O$.
Let $v\in V[\la_v], w\in W[\la_w]$. Assume that $v, w$ are homogeneous with respect to the grading.
Consider
$$
M^+_\la  \buildrel \Phi^v_\la \over \longrightarrow M^+_{\la-\la_v}\T V
\buildrel \Phi^w_{\la -\la_v}\T\,1  \over 
\longrightarrow M^+_{\la-\la_v-\la_w}\T W\T V
$$
and denote this composition $\Phi^{w,v}_\la$. 

Define the main object of this paper, a linear operator $J_{W,V}(\la):W\T V \to W\T V$
as follows. Find $u\in W\T V [\la_v+\la_w]$ such that $\Phi_\la^{w,v}=\Phi_\la^u$
and set 
\begin{equation}\label{dj}
J_{W,V}(\la) \,w\T v = u.
\end{equation}

\begin{lemma}
$J_{W,V}(\la)$ is strictly upper triangular, i.e. $J_{W,V}(\la) \,w\T v = w\T v +
\sum w_i\T v_i$ where deg $w_i\,<\,$ deg $w$. $\square $
\end{lemma}
\begin{corollary}  $J_{W,V}(\la) = 1 + N$,
where $N$ is locally nilpotent,  hence $J_{W,V}(\la)$ is invertible
and $J_{W,V}^{-1}(\la) = 1 - N + N^2 -....$

\end{corollary}

We call the operators $J_{W,V}(\la)$ {\it fusion matrices}.

Define {\it a quantum dynamical R-matrix} 
$R_{V,W}(\la): V\T W \to V\T W$ by
\begin{equation}\label{r}
R_{V,W}(\la)=J_{V,W}^{-1}(\la)\Rm^{21}|_{V\T W} J^{21}_{W,V}(\la).
\end{equation}
\begin{thm} Let $v\in V, w\in W$ be homogeneous elements with respect
to the grading and $A_0$. Let $R_{V,W}(\la)\, v\T w= \sum_i v_i\T w_i$
where $v_i, w_i$ are homogeneous too. Then
\begin{equation}\label{}
(1\,\T P \Rm|_{W\T V})\, \Phi_\la^{w,v} =
\sum_i \, \Phi_\la^{v_i,w_i}
\end{equation}
where $P$ is the operator of permutation.
\end{thm}
The proof follows from the definition of $R_{V,W}(\la)$.
\begin{thm}\label   {R}

\begin{enumerate} 
\item[I.] $J$ satisfies the  2-cocycle condition,
\begin{equation}\label{2c}
J_{V\T W, U}(\la)\,(J_{V, W}(\la - h^{(3)})\T\,1)=
J_{V,  W\T U}(\la)\,(1\,\T J_{W, U}(\la)).
\end{equation}
\item[II.] $R (\la)$ satisfies the quantum dynamical Yang-Baxter equation (QDYB),
\bean\label{qdyb}
R^{12}(\la -  h^{(3)})\,
R^{13}(\la )\,
R^{23}(\la -  h^{(1)})\,
=
R^{23}(\la )\,
R^{13}(\la -  h^{(2)})\,
R^{12}(\la)\,.
\eean
\end{enumerate}
\end{thm}
In these formulae $R^{12}(\la -  h^{(3)})\,v\T w\T u=
(R(\la - \la_u)\,v\T w)\,\T u$ if $u\in U[\la_u]$ and other symbols are
defined analogously.

\begin{proof}
The first statement of the Theorem is trivial. To prove the second statement
define 
$$
\Phi_\la^{v_n,...,v_1}:M^+_\la \to M^+_{\la-\sum_{i=1}^n\la_{v_i}}\T V_n \T ...
\T V_1
$$
as the composition
\bean
M^+_\la \buildrel \Phi^{v_1}_\la \over \longrightarrow M^+_{\la-\la_{v_1}}\T V_1
\buildrel \Phi^{v_2}_{\la -\la_{v_1}}\T\,1  \over 
\longrightarrow M^+_{\la-\la_{v_1}-\la_{v_2}}\T V_2\T V_1 \to ...
\notag
\eean
In other words, 
$$
\Phi_\la^{v_n,...,v_1}= (\Phi^{v_n}_{\la - \sum_{i=1}^{n-1}\la_{v_i}}\T 
\, 1^{n-1}) \cdots (\Phi^{v_2}_{\la - \la_{v_1}}\T \, 1)\,\Phi^{v_1}_{\la}.
$$
\begin{lemma}
$$
\Phi_\la^{v_n,...,v_1}=\Phi_\la^{v_n,...,v_{i+2}, J_{V_{i+1},V_i}(\la-\sum_{j=1}^{i-1}
\la_{v_j})\, v_{i+1}\T v_i, v_{i-1},...,v_1}
$$
\end{lemma}
The proof is by definition of $J$.

Let $a\in V, b\in W, c \in U$.
Then 
$$
\Phi_\la^{a,b,c}=\Phi_\la^{a,J_{W,U}(\la)\,b\,\T\, c}=
\Phi_\la^{ J_{V, W\T U}(\la)\, (1\, \T \,J_{W,U}(\la))\, a\,\T \,b\,\T\, c},
$$
$$
\Phi_\la^{a,b,c}=\Phi_\la^{J_{V,W}(\la - \la_c)\,a\,\T \, b,  c}=
\Phi_\la^{ J_{V\T W, U}(\la)\, (J_{V,W}(\la - \la_c)\,\T\,1)\, a\,\T \,b\,\T\, c}.
$$
This proves the first statement of the Theorem.

For $y\in V_n\T ... \T V_1 [\nu],\,{}\, y=\sum_i
v^i_n\T ... \T v^i_1$, set
$\Psi^y_{V_n,...,V_1}(\la)= \sum_i \Phi_\la^{v^i_n, ... v^i_1}.$
\begin{lemma}
$$
P_{V_{i+1}, V_i}\Rm_{V_{i+1}, V_i} \Psi^{v_n\,\T\, ... \,\T\, v_1}_{V_n,...,V_1}(\la)=
\Psi^{R_{V_i,V_{i+1}}(\la-\sum_{j=1}^{i-1}\la_{v_j})
P_{V_{i+1},V_i}\,v_n\,\T\, ... \,\T\, v_1}
_{V_n,...,V_{i},V_{i+1},...,V_1}(\la).
$$
\end{lemma}
The proof is by definition of the quantum dynamical R-matrix.

In order to prove the second statement of the Theorem
we apply the Lemma to the case $n=3$. Namely, we consider the function 
$\Psi_{V_3,V_2,V_1}^{v_3\T v_2\T v_1}(\la)$ and express it 
via $\Psi_{V_1,V_2,V_3}^{w_1\T w_2\T w_3}(\la)$ in two different ways, 
using the two different reduced decompositions of the permutation
$123\to 321$. Comparing the two answers, we get the 
Theorem. 

\end{proof}

{\bf  Remark.} The explicit form of $J_{VW}$ has been recently 
computed in [A]. The fact that the twist in [A] coincides with our $J_{VW}$
is proved in Chapter 9. 

%\section{The Tensor Functor }\label{tf}

\subsection{The tensor functor and exchange matrices}\label{cvs}
We recall (in a slightly generalized form) the setting of Chapter 3 of 
\cite{EV2}. 
Let $A_0$ be a commutative, cocommutative finitely generated Hopf algebra
such that the group $T=$Spec $A_0$ is connected.  
Introduce a category $\V$
of $A_0$-vector spaces as follows. 

The objects of $\V$ are diagonalizable $A_0$ modules,
$V=\oplus_{\la \in T}V[\la], \, V[\la]=\{v\in V\,|\, a_0v=\la(a_0)v\}.$

Let $M_T$ be the field of meromorphic functions on $T$ and $V, W \in \V$.
Define the space $\Hom_\V(V, W)$ as the space 
$\Hom_{A_0}(V, M_T\T  W)$, thus a homomorphism of $V$ to $W$ (for finite
dimensional  $V,W \in \V$) is
a meromorphic function on $T$ with values in $\Hom_{A_0}(V, W)$.

Define a tensor structure on the category $\V$. Namely, let the tensor product
of two objects $V\TB W$ be the standard tensor product of two diagonalizable $A_0$ modules.
Define the tensor product $\bar {\T}$ of two morphisms $f: V \to V'$ and $g:W\to W'$ as
\begin{equation}\label{bar}
f\bar{\T} g (\la)= f^{(1)}(\la-h^{(2)})(1\,\T g(\la))
\end{equation}
where $f^{(1)}(\la-h^{(2)})(1\,\T g(\la)) \, u\T v=
(f(\la-\mu) u) \,\T\,g(\la) v $ if $g(\la) v \in W'[\mu]$.

Let $A$ be a nondegenerate polarized Hopf algebra.
Consider the category $\cal O$ of graded $A$-modules bounded from above and
diagonalizable over $A_0$. We construct a tensor functor from the category $\Oc$
to the category $\V$. 

By definition a tensor functor from $\Oc $ to $\V$ is a functor $F:\Oc\to \V$
and for any $V,W \in \Oc$ an isomorphism $J_{V,W}: F(V)\TB F(W) \to F(V\T W)$ 
such that $\{J_{V,W}\}$ is functorial and the two compositions
$F(U)\TB F(V)\TB F(W)\to F(U\T V)\TB F(W)\to F(U\T V\T W)$ and
$F(U)\TB F(V)\TB F(W)\to F(U)\TB F(V\T W)\to F(U\T V\T W)$ coincide.
Then $J$ is called {\it a tensor structure on} $F$.

Define a tensor functor $F : \cal O \to \V$ by sending
an object $V\in \cal O$ to $F(V)= V$, considered as an $A_0$-module,
and sending an $A$-homomorphism $\al : V\to W$ to $F(\al)=\al : V\to W$.

Define a tensor structure on $F$ by
\begin{equation}\label{ts}
J_{V,W}(\la): F(V)\TB F(W) \to F(V\T W)
\end{equation}
where $J_{V,W}(\la)$ is defined by \Ref{dj}.
\begin{lemma}
Formula \Ref{ts} defines a tensor structure on $F$, i.e. the two compositions
$F(U)\TB F(V)\TB F(W)\to F(U\T V)\TB F(W)\to F(U\T V\T W)$ and
$F(U)\TB F(V)\TB F(W)\to F(U)\TB F(V\T W)\to F(U\T V\T W)$ coincide.
\end{lemma}
\begin{proof} The statement of the Lemma is equivalent to formula
\Ref{2c}.
\end{proof}

Define a braiding in $\cal O$ by $\bt=P \Rm$. Introduce
\begin{equation}\label{br}
F(\bt): F(V)\TB F(W)\to F(W)\TB F(V)
\notag
\end{equation}
as the composition
\bean
F(V)\TB F(W) \buildrel J_{V,W}(\la) \over \longrightarrow 
F(V\T W) \buildrel P\Rm_{V,W}  \over 
\longrightarrow F(W\T V) 
\buildrel J^{-1}_{W,V}(\la) \over \longrightarrow 
F(W)\TB F(V).
\eean
Thus we have
$F(\bt)_{V,W}=J^{-1}_{W,V}(\la)P_{V,W}\Rm |_{V\T W}J_{V,W}(\la)$. 
In particular,
$F(\bt)_{V,W}P_{W,V}=J^{-1}_{W,V}(\la)\Rm^{21} |_{V\T W}J^{21}_{V,W}(\la)
=R_{W,V}(\la)$, cf. \Ref{r}.
Notice that in Theorem \ref{R} we showed that the R-matrix $R(\la)$ satisfies the QDYB
equation, now it also follows from this tensor category construction and Theorem 3.3
in [EV2].

The operators $R_{V,W}(\la)$ will be called {\it the exchange matrices}.

\section{ $H$-Hopf algebroids}

\subsection{ Definitions }\label{dfns} 
In this Section we recall the definition
of an $H$-Hopf algebroid, cf. [EV2].
Let $H$ be a commutative and cocommutative Hopf algebra over $\C$,
$T=$Spec $H$   a commutative affine algebraic group.
%$\ti$ the Lie algebra of $T$. Assume that $T$ is connected. 
Let $M_{T}$ denote the field of meromorphic functions on $T$.
{\it An $H$-algebra } is an associative algebra
$A$ over $\C$ with $1$, endowed with an $T$-bigrading
$A=\oplus_{\al,\beta\in T}A_{\al\beta}$
(called the weight decomposition), and two algebra embeddings
$\mu_l,\mu_r:M_{T}\to A_{00}$ (the left and the right moment maps), such
that for any $a\in A_{\al\beta}$ and $f\in M_{T}$, we have
\begin{equation}\label{mmm}
\mu_l(f(\la))a=a\mu_l(f(\la+\al)),\quad
\mu_r(f(\la))a=a\mu_r(f(\la+\beta)).
\end{equation}
Here $0\in T$ denotes the unit element and $\la + \al$ denotes the sum  in $T$.

{\it A morphism} $\phi:A\to B$ of two $H$-algebras
is an algebra homomorphism, preserving
the moment maps. 

Let $A,B$ be two $H$-algebras
and $\mu_l^A,\mu_r^A,\mu_l^B,\mu_r^B$ their moment maps.
Define their {\it matrix tensor product}, $A\wo B$, which is also an $H$-algebra.
Let
\begin{equation}\label{pr}
(A\wo B)_{\al\delta}:=\oplus_{\beta}A_{\al\beta}
\T_{M_{T}} B_{\beta\delta},
\end{equation}
 where $\T_{M_{T}}$ means
the usual tensor product
modulo the relation $\mu_r^A(f)a\T b=a\T \mu_l^B(f)b$, for any
$a\in A,b\in B, f\in M_{T}$.
Introduce a multiplication in $A\wo B$ by the rule $(a\T b)(a'\T b')=
aa'\T bb'$.  Define
the moment maps for $A\wo B$ by
$\mu_l^{A\wo B}(f)=\mu_l^A(f)\T 1$,
$\mu_r^{A\wo B}(f)=1\T \mu_r^B(f)$.

{\it A coproduct} on an $H$-algebra $A$
is a homomorphism of $H$-algebras
$\Delta: A\to A\wo A$.

Let $D_T$ be the algebra of difference operators
$M_{T}\to M_{T}$, i.e.
the operators of the form $\sum_{i=1}^nf_i(\la)\Ti_{\beta_i}$, where
$f_i\in M_{T}$, and
for $\beta\in T$ we denote by $\Ti_\beta$ the field automorphism
of $M_{T}$ given by $(\Ti_\beta f)(\la)=f(\la+\beta)$.

The algebra $D_T$ is an  example of an $H$-algebra  
if we define the weight decomposition by $D_T=\oplus (D_T)_{\al\beta}$,
where $(D_T)_{\al\beta}=0$ if $\al\ne\beta$, and
$(D_T)_{\al\al}=\{f(\la)\Ti_\al^{-1}: f\in M_{T}\}$, and the moment maps
$\mu_l=\mu_r:M_{T}\to (D_T)_{00}$
to be the tautological isomorphism.

For any $H$-algebra $A$, the algebras
$A\wo D_T$ and $D_T \wo A$ are canonically isomorphic to $A$.
In particular, $D_T$ is canonically isomorphic to $D_T\wo D_T$. 
Thus the category of $H$-algebras
equipped with the product $\wo$ is a
monoidal category, where the unit object is $D_T$.

{\it A counit } on an $H$-algebra $A$ is
a homomorphism of $H$-algebras $\epe: A\to  D_T$.

{\it An $H$-bialgebroid } is a
$H$-algebra $A$ equipped with a coassociative coproduct
$\Delta$ (i.e. such that $(\Delta\T \Id_A)\circ \Delta=
(\Id_A\T \Delta)\circ \Delta$, and a counit
$\epe$ such that $(\epe\T \Id_A)\circ \Delta=
(\Id_A\T \epe)\circ \Delta=\Id_A$.

For example, $D_T$ is an $H$-bialgebroid where $\Dl : D_T \to D_T\Tb D_T$
is the canonical isomorphism and $\epe = \Id$.

Let $A$ be an $H$-algebra. A linear map $S:A\to A$ is called {\it an antiautomorphism}
of an $H$-algebra if it is an antiautomorphism
of algebras and $\mu_r\circ S=\mu_l,\, \mu_l\circ S=\mu_r$. From 
these conditions it follows that $S(A_{\al\beta})=A_{-\beta,-\alpha}$.

Let $A$ be an $H$-bialgebroid, and let $\Delta$, $\epe$
be the coproduct and counit of $A$. For $a\in A$, let
\begin{equation}\label{pres'}
\Delta(a)=\sum_i a^1_i\T a^2_i.
\end{equation}

{\it An antipode} on the $H$-bialgebroid $A$
is an antiautomorphism of $H$-algebras $S:A\to A$ such that
for any $a\in A$ and any presentation \Ref{pres'} of $\Delta(a)$, one has
\begin{equation}
\sum_i a_i^1S(a_i^2)=\mu_l(\epe(a)1),\
\sum_i S(a_i^1)a_i^2=\mu_r(\epe(a)1),
\notag
\end{equation}
where $\epe(a)1\in M_{T}$ is the result of the application of the difference
operator $\epe(a)$ to the constant function $1$.

An $H$-bialgebroid with an antipode is called {\it an $H$-Hopf algebroid.}

Let $W$ be a diagonalizable $H$-module, 
$W=\oplus_{\la \in T} W[\la], \, W[\la]=\{w\in W\,|\, aw=\la(a)w, \text{for all}\,
a\in H \}$,  and let $D^\al_{T, W}\subset
\text{Hom}_\C(W,W\T D_T)$ be the space of all difference
operators on $T$ with coefficients in $\End_\C(W)$, which have weight $\al \in T$ with
respect to the action of $H$ in $W$.

Consider the algebra $D_{T, W}=\oplus_\al D_{T, W}^\al$.
This algebra has a weight decomposition
$D_{T,W}=\oplus_{\al,\beta} (D_{T,W})_{\al\beta}$ defined as follows:
if $g\in \text{Hom}_\C(W,W\T M_{T})$ is
an operator of weight $\beta-\al$,
then $g\Ti_{\beta}^{-1}\in (D_{T,W})_{\al\beta}$.

Define the moment maps $\mu_l,\mu_r: M_{T}\to (D_{T, W})_{00}$
by the formulas $\mu_r(f(\la))=f(\la)$,
$\mu_l(f(\la))=f(\la- h)$ where $f(\la - h)w=f(\la - \mu)w$
if $w\in W[\mu],\, \mu \in T$.  The algebra $D_{T, W}$ equipped with this weight
decomposition and these moment maps is an $H$-algebra.

Let $f\in \text{Hom}(W,W\T M_{T})$ and $g\in \text{Hom}(U, U\T M_{T})$.
Define $f\bo g \in \Hom (W\T U, W\T U \T M_{T})$ as
\begin{equation}\label{bar'}
f\bar{\T} g (\la)= f^{(1)}(\la-h^{(2)})(1\,\T g(\la))
\end{equation}
where $f^{(1)}(\la-h^{(2)})(1\,\T g(\la)) \, w\T u=
(f(\la-\mu) w )\,\T\,g(\la) u $ if $g(\la) u \in U[\mu]$, cf. \Ref{bar}.

\begin{lemma} $\text {[EV2]}$  There is a natural embedding of $H$-algebras
$\theta_{W,U}:D_{T, W}\wo D_{T, U}\to D_{T, W\T U}$,
given by the formula $f\Ti_\beta\,\T g\,\Ti_\delta\to (f\bo g)\Ti_\delta$.
This embedding is  an isomorphism if $W,U$ are finite-dimensional.
\end{lemma}

{\it A dynamical representation} of an
$H$-algebra $A$ is a diagonalizable $H$-module
$W$ endowed with a homomorphism of $H$-algebras
$\pi_W:A\to D_{T, W}$.
{\it A homomorphism} of dynamical representations
$\phi: W_1\to W_2$ is an element of $\text{Hom}_{\C}(W_1,W_2\T M_{T})$
such that $\phi\circ \pi_{W_1}(x)=
\pi_{W_2}(x)\circ \phi$ for all $x\in A$.

{\bf Example.} If $A$ has a counit, then
$A$ has {\it the trivial representation}: $W=\C$, $\pi=\epe$.

If $A$ is an $H$-bialgebroid, $W$ and $U$ are two dynamical representations of $A$,
then the $H$-module
$W\T U$ is a dynamical representation,
 $\pi_{W\T U}(x)=\theta_{WU}\circ (\pi_W\T \pi_U)\circ \Delta(x)$.
If $f:W_1\to W_2$ and $g:U_1\to U_2$
are homomorphisms of dynamical representations, then so is
$f\bo g :W_1\T U_1\to W_2\T U_2$.
Thus, dynamical representations of $A$ form a monoidal category
$\text{Rep}(A)$,
whose identity object is the trivial representation.

%The category $\text{Rep}(A)$ is equipped with a natural
%tensor functor $\text{Rep} (A)\to $ to the category of
%$\h$-vector spaces -- the forgetful functor.

If $A$ is an $H$-Hopf algebroid and $V$ is a dynamical representation,
then one can define the left and right dual dynamical representations
${}^*W$ and $W^*$ as follows, see [EV2]. 

If $(W, \pi_W)$ is a dynamical representation of an $H$-algebra $A$,
we denote  $\pi_W^0:A\to \text{Hom}(W,W\T M_T)$
the map defined by $\pi_W^0(x)w=
\pi_W(x)w$, $w\in W$
 (the difference operator $\pi_W(x)$ restricted to the constant functions).
It is clear that $\pi_W$ is completely determined by $\pi_W^0$.

Let $(W,\pi_W)$ be a dynamical representation of $A$.
Then {\it the right dual representation} to $W$ is $(W^*,\pi_{W^*})$,
where $W^*$ is the $H$-graded dual to $W$, and
\bean\label{rdual}
\pi_{W^*}^0(x)(\la)=\pi_W^0(S(x))(\la+ h- \alpha)^t
\eean
for $x\in A_{\al\beta}$, where $t$ denotes dualization.
The left dual representation to $W$ is $({ }^*W,\pi_{{ }^*W})$,
where ${ }^*W=W^*$, and
\bean\label{ldual}
\pi_{{ }^*W}^0(x)(\la)=
\pi_W^0(S^{-1}(x))(\la+ h - \alpha)^t
\eean
for $x\in A_{\al\beta}$.
Here $(S(x))(\la+ h- \alpha)^t$ denotes
the result of two operations applied successively to $S(x)$: 
shifting of the argument, and dualization. Similarly, 
$(S^{-1}(x))(\la+ h - \alpha)^t$ denotes the result of
of three operations applied successively to $S(x)$: inversion,
shifting of the argument, and dualization.

Formulas  \Ref{rdual} and \Ref{ldual} define dynamical representations of $A$.
Moreover, if $A(\la):W_1\to W_2$ is a morphism
of dynamical representations, then $A^*(\la):=A(\la+ h)^t$
defines  morphisms $W_2^*\to W_1^*$ and ${ }^*W_2\to { }^*W_1$.

\subsection{ An  $H$-bialgebroid associated to a
 function $R: T \to \End(V\T V)$}\label{hh}
In this Section we recall a construction from [EV2]
of an $H$-bialgebroid $A_R$ associated to a meromorphic
function $R: T \to \End(V\T V)$ where $V$ is a finite dimensional
diagonalizable
$H$-module and $R(\la)$ is invertible for generic  $\la$.

By definition the algebra $A_R$  is generated by two copies of
$M_{T}$ (embedded as subalgebras) and  matrix elements of the operators
$L^{\pm 1}\in \text{End}(V)\T A_R$.
 We denote the elements of the first copy of
$M_{T}$ by $f(\la^1)$ and of the second copy by $f(\la^2)$,
where $f\in M_{T}$. We denote  $(L^{\pm 1})_{\al\beta}$
the weight components of $L^{\pm 1}$ with respect to the natural
$T$-bigrading on $\End (V)$, so
that $L^{\pm 1}=(L^{\pm 1}_{\al\beta})$, where
$L^{\pm 1}_{\al\beta}\in \text{Hom}_\C(V[\beta],V[\al])\T A_R$.

Introduce the moment maps for $A_R$ by
$\mu_l(f)=f(\la^1)$, $\mu_r(f)=f(\la^2)$, and the weight
decomposition by $f(\la^1),f(\la^2)\in (A_R)_{00}$, $L_{\al\beta}\in
\text{Hom}_\C(V[\beta],V[\al])\T (A_R)_{\al\beta}$.

The defining relations for $A_R$ are:
\bean
f(\la^1)L_{\al\beta}=L_{\al\beta}f(\la^1+\al);\
\qquad
f(\la^2)L_{\al\beta}=L_{\al\beta}f(\la^2+\beta);
\\
LL^{-1}=L^{-1}L=1;
\qquad
[f(\la^1),g(\la^2)]=0;
\eean
and the dynamical Yang-Baxter relation
\bean\label{ybr}
R^{12}(\la^1)L^{13}L^{23}=:L^{23}L^{13}R^{12}(\la^2):.
\eean
Here the :: sign 
means that the matrix elements of $L$
should be put on the right of the matrix elements of $R$. Thus, if
$\{v_a\}$ is a homogeneous basis of $V$, and $L=\sum E_{ab}\T L_{ab}$,
$R(\la)(v_a\T v_b)=\sum R^{ab}_{cd}(\la)v_c\T v_d$, then \Ref{ybr} has the form
$$
\sum R^{xy}_{ac}(\la^1)L_{xb}L_{yd}=\sum R_{xy}^{bd}(\la^2)
L_{cy}L_{ax},
$$
where we sum over repeated indices.

Define the coproduct on $A_R$, $\Delta: A_R\to A_R\Tb A_R$, by
$$
\Delta(L)=L^{12}L^{13},\Delta(L^{-1})=(L^{-1})^{13}(L^{-1})^{12}.
$$
Define the counit by 
$$
\epe(L_{\al\beta})=
\delta_{\al\beta} \Id_{V[\al]} \T \Ti_\al^{-1},
\epe((L^{-1})_{\al\beta})=
\delta_{\al\beta} \Id_{V[\al]} \T \Ti_\al,
$$
where $\Id_{V[\al]}: V[\al]\to V[\al]$ is the identity operator.
On an antipode for $A_R$ see Section 4.5 in [EV2].

{\bf A rational  $H$-bialgebroid associated to a
rational function $R: T \to \End(V\T V)$.}
Assume that a function $R: T \to \End(V\T V)$ is a rational
function of $\la$,
where $V$ is a finite dimensional diagonalizable
$H$-module and $R(\la)$ is invertible for generic  $\la$.

The  $H$-bialgebroid $A_R$ is defined  over
the field of meromorphic functions $M_T$. We replace the field of meromorphic
functions $M_T$ by the field of rational functions $\C(T)$
and define in the same way {\it the rational $H$-bialgebroid } $A_{rat,R}$ 
associated to a rational function $R$.

\section{ The Exchange Dynamical Quantum Groups}

\subsection{The definition of an exchange dynamical quantum group}
Let $A$ be a polarized and nondegenerate Hopf algebra as in Section 1.
Assume that $T=$Spec $A_0$ is connected.  

%Let $\ti$ be the Lie algebra of $T$. Consider on $\h:=\ti^*$
%the abelian Lie algebra structure, $[\al,\bt]=0$ for all $\al, \bt \in \h$.

Let $\Rm\in A\,\hat {\T}\, A$ be  a quasitriangular structure  on $A$. We always assume
that $\Rm\in A_{\geq 0}\,\hat {\T}\, A_{\leq 0}$.

Let $\Oc_0\subset \Oc$ be a full abelian tensor subcategory which is semisimple
and such that all modules in $\Oc_0$ are finite dimensional.
(Remind that a full subcategory $\Oc_0$ consists of some objects of $\Oc$ and for any
$V, W \in \Oc_0$ we have $\Hom_{\Oc_0}(V,W)=\Hom_{\Oc}(V,W)$.) 
Let $Ir \subset \Oc_0$ be the set of all irreducible modules.

%Let $U \in \Oc_0$ and $
%Let exp : $\ti \to T$ be the exponential map and 

Examples of such categories $\Oc_0$ are provided by the categories of
finite dimensional representations of semisimple Lie algebras and corresponding
quantum groups ( not at roots of unity).

The goal of this Section is to define an $A_0$-Hopf algebroid $E=E(\Oc_0)$
called {\it an exchange dynamical quantum group}.

Define $E$ as a vector space to be
$$
\Mm \T_\C \Mm \T_\C \E
$$
where $\E=\oplus_{U\in Ir}U\T U^*$ and $U^*$ is the  dual module to $U$. 
A $T$-bigrading on $E$ is defined by
$E=\oplus_{\al,\bt \in T}E_{\al,\bt}$,
where $E_{\al,\bt}=\Mm \T_\C \Mm \T_\C \E_{\al,\bt}$ and
$\E_{\al,\bt} \subset \oplus_{U\in Ir}U\T U^*$ is the subspace
generated by all elements of the form $u \T v \in U[\al] \T (U [\bt])^*,\,U\in Ir$.

Let $\C_A$ be the trivial $A$-module, $\C_A=\C e $.
The subspace $E_{0,0}$ has a component coming from the trivial module,
$ \Mm \T \Mm \T \C_A \T \C_A^*$. For a meromorphic function $f(\la)\in \Mm$,
the elements $f(\la) \T 1 \T e \T e^*$ and $1\T f(\la) \T e \T e^*$ will be denoted
$f(\la^1)$ and $f(\la^2)$, respectively.

Let $v_i^U$ be a basis in $U \in Ir$, which is homogeneous with respect
to $T$ and the $\Z$-grading. Then $v_i^U \T (v_j^U)^*$ form a basis in $\E$.
Let $\om^U_{i}\in T$ be the weight of $v^U_{i}$.

Set $L^U_{ij}=1\T 1\T v_i^U \T (v_j^U)^*$. Define a linear map
$E^U_{ij}:U \to U$ by $E^U_{ij} v^U_{k}=\dl_{jk}v^U_i$. Introduce 
$L^U\in \End\,(U)\T \E$ by
$$
L^U=\sum_{ij} E^U_{ij}\T L^U_{ij}.
$$
The relations in $E$ between $f(\la^1), f(\la^2)$, and $L^U_{ij}$
are defined by
\bean \label{ffr} 
f(\la^1)f(\la^2)=f(\la^2)f(\la^1), 
\\
\label{flr}
f(\la^1)L^U_{ij}=L^U_{ij}f(\la^1+\om^U_{i}),
\qquad
f(\la^2)L^U_{ij}=L^U_{ij}f(\la^2+\om^U_{j})
\eean
In order to define the product of two elements  $L^V_{ij}$ and
$L^W_{i'j'}$ we will consider
$(L^V)^{23}(L^W)^{13}\in \End\, (V)\T \End\,(W)\T E$.

Let $U\in Ir, \, V, W
\in \Oc_0$ and $H^U_{V,W}=\Hom_A (U,V\T W)$.
Then we have an isomorphism $\tau_{V,W}:\oplus_{U\in Ir} H^U_{V,W}\T U \to V\T W$
given by $\tau_{V,W} (h\T u)= h(u)$. Let $\ta ^U_{V,W}:V\T W \to H^U_{V,W}\T U$
be the projection along the other summands,
$\tau^U_{V,W}:  H^U_{V,W}\T U \to V\T W$
the restriction of $\tau_{V,W}$ to the isotypic component
$H^U_{V,W}\T U$. We have
$$
\ta ^U_{V,W} \tau ^U_{V,W} = \Id, \qquad
\tau ^U_{V,W} \ta ^U_{V,W} = p_U
$$
where $p_U$ is the projection on the $U$-isotypical component.

Define the product of elements in the exchange quantum group by a formula
analagous to the formula for the product of matrix elements of representations of a
group considered as functions on the group. Namely, define the product
$L^V_{ij}L^W_{i'j'}$ by
\bean\label{dr}
(L^V)^{23}(L^W)^{13}\,=\,: (J_{W,V}^{12}(\la^1))^{-1}\,
\sum_{U \in Ir} (\tau^U_{W,V})^{12}\,(\Id_{H^U_{W,V}}\T L^U)\,
(\ta _{W,V}^U)^{12}\,J^{12}_{W,V}(\la^2):\,.
\eean
This is an identity in $\End\, (W)\T \End\,(V)\T E$.
Here the :: sign (``normal ordering'')
means that the matrix elements of 
$L^U$ should be put on the right of the matrix elements of $J_{V,W}(\la^1),
 J_{V,W}(\la^2)$. Thus, if
$(J_{V,W}^{12}(\la^1))^{-1}= E^V_{ij}\T E^W_{kl} \T a_{ijkl}(\la^1)$,
$\sum_{U \in Ir} (\tau^U_{V,W})^{12}\,(\Id_{H^U_{V,W}}\T L^U)\,
(\ta _{V,W}^U)^{12}=E^V_{i"j"}\T E^W_{k"l"} \T a"_{i"j"k"l"}$,
$J_{V,W}^{12}(\la^2)= E^V_{i'j'}\T E^W_{k'l'} \T a'_{i'j'k'l'}(\la^2)$,
then \Ref{dr} has the form
\bean
L^W_{kl'} L^V_{ij'}=\sum a_{ijkl}(\la^1)
a'_{ j"j'l"l'}(\la^2)a"_{jj"ll"}.
\notag
\eean
More generally, let $a=a_1...a_n$ be a monomial in generators of $E$;
so each of the factors has the form $f(\la^1), f(\la^2)$, or $L^V_{ij}$.
Define the normal ordering $:a:$ as the product of the same elements
$a_1,...,a_n$ in which all elements of the form
$f(\la^1), f(\la^2)$ are put on the left and the remaining elements
of the form $L^V_{ij}$ are put on the right in the same order as in $a$.
Extend by linearity the normal ordering operation to all polynomials
in generators in $E$. If $v=(v_1,...,v_l)$ is a vector whose
coefficients are polynomials in generators of $E$, then
define  the normal ordering $:v:$ as $:v:=(:v_1:,...,:v_l:)$.

Let $\C_A$ be the trivial module. Since $J_{\C_A, V}=J_{V,\C_A}=\Id_V$ we have
\bean
(L^{\C_A})^{23}(L^W)^{13}\,=\,(L^W)^{13}, \qquad
(L^{V})^{23}(L^{\C_A})^{13}\,=\,(L^V)^{23}.
\eean
\begin{corollary}
The element $1\T 1\T e\T e^*$ of $E$ corresponding  to the trivial 
module is the unit element of the algebra $E$.
\end{corollary}

\begin{thm}
$E$ is an associative algebra.
\end{thm}
\begin{proof} We start with preliminary lemmas.
\begin{lemma}\label{funct}

\begin{enumerate}
\item[I.] Let $V'\subset V$ be objects in $\Oc_0$, then
$J_{V,W}|_{V'\T W}=J_{V',W}$.
 Let $W'\subset W$ be $A$-modules, then
$J_{V,W}|_{V\T W'}=J_{V,W'}$.
\item[II.]  Let $V=V_1\oplus V_2$, then
$J_{V,W}=J_{V_1,W}\oplus J_{V_2,W}.$
 Let $W=W_1\oplus W_2$, then
$J_{V,W}=J_{V,W_1}\oplus J_{V,W_2}.$
\item[III.]
For $U\in Ir, \,V,W,Z \in \Oc_0$,
the maps $Z\T W\T V \to H^U_{Z,W}\T U\T V$ given by
$(\ta ^U_{Z,W}\T \Id_V)\,J_{Z\T W, V}$ and
$(\Id_{H^U_{Z,W}}\T J_{U,V})\,(\ta ^U_{Z,W} \T \Id_V)$
coincide.
\item[IV.] The maps $Z\T H^U_{W,V}\T U \to Z\T W\T V$
given by $J_{Z, W\T V}\,(\Id_Z \T \tau^U_{W,V})$
and $(\Id_Z \T \tau^U_{W,V})\,(J_{Z,U})^{13}$
coincide. In particular, 
$J^{-1}_{Z, W\T V}\,(\Id_Z \T \tau^U_{W,V})=
(\Id_Z \T \tau^U_{W,V})(J^{-1}_{Z,U})^{13}$
\end{enumerate}
\end{lemma}
The Lemma follows from functorial properties of $J$.
%\begin{corollary}
%$(\tau^U_{Z,W}\T \Id_V)\,(\Id_{H^U_{Z,W}}\T J^{-1}_{U,V})=(J_{Z\T W, V})^{-1}
%(\tau^U_{Z,W}\T \Id_V).$
%\end{corollary}
Now we prove the Theorem. We want to show that
\bean\label{!}
(L^V)^{34}((L^W)^{24}(L^Z)^{14})=
((L^V)^{34}(L^W)^{24})(L^Z)^{14}.
\eean
We have
\bean\label{!!}
\text{RHS}=:(J^{-1}_{W,V}(\la^1))^{23}\sum_{U\in Ir}(\tau^U_{W,V})^{23}
(\Id_{H^U_{W,V}}\T L^U)^{234}(\ta ^U_{W,V})^{23} J^{23}_{W,V}(\la^2)
(L^Z)^{14}:\,.
\eean
First we replace $(\ta ^U_{W,V})^{23} J^{23}_{W,V}(\la^2)(L^Z)^{14}$
with $(L^Z)^{14}(\ta ^U_{W,V})^{23} J^{23}_{W,V}(\la^2)$.
Consider $(\Id_{H^U_{W,V}}\T L^U)^{234}(L^Z)^{14}$ as an element
of the tensor product $\End (Z)\T \End (H^U_{W,V})\T \End (U)\T \E$,
then the element  $(\Id_{H^U_{W,V}}\T L^U)^{234}(L^Z)^{14}$
takes the form $(\Id_{H^U_{W,V}})^{2'} (L^U)^{3'4}(L^Z)^{14}$
where ${}^{2'}, {}^{3'}$ label these new tensor factors.
Applying  formula \Ref{dr} to the first, third and fourth factors, we get
\bean
  (L^U)^{3'4}(L^Z)^{14} 
\,=\,
\notag
\\
: (J_{Z,U}^{13'}(\la^1))^{-1}\,
\sum_{Y \in Ir} (\tau^Y_{Z,U})^{13'}\,(\Id_{H^Y_{Z,U}}\T L^Y)\,
(\ta _{Z,U}^Y)^{13'}\,J^{13'}_{Z,U}(\la^2):\,.
\notag
\eean
Returning to \Ref{!!} we get
\bean
\text{RHS}=:(J^{-1}_{W,V}(\la^1))^{23}\sum_{U\in Ir}(\tau^U_{W,V})^{23}
\notag
\\
 (J^{13}_{Z,U}(\la^1))^{-1}
\sum_{Y \in Ir} (\tau^Y_{Z,U})^{13}\,(\Id_{H^Y_{Z,U}}\T L^Y)\,
(\ta _{Z,U}^Y)^{13}\,J^{13}_{Z,U}(\la^2)
(\ta ^U_{W,V})^{23} J^{23}_{W,V}(\la^2): \,.
\notag
\eean
Applying Lemma \ref{funct} we get
\bean\label{RHS}
\text{RHS}=:(J^{-1}_{W,V}(\la^1))^{23}(J^{-1}_{Z,W\T V}(\la^1))^{1,23}
\sum_{U\in Ir}(\tau^U_{W,V})^{23}
\sum_{Y \in Ir} (\tau^Y_{Z,U})^{13}\,(\Id_{H^Y_{Z,U}}\T L^Y)\,\times
\notag
\\
(\ta _{Z,U}^Y)^{13}\,
(\ta ^U_{W,V})^{23} 
J^{1,23}_{Z,W\T V}(\la^2)
J^{23}_{W,V}(\la^2): \,.
\eean
Now we compute the left hand side of \Ref {!}.
\bean\label{LHS}
\text{LHS}=\sum_{ij}(E^V_{ij})^3(L^V_{ij})^4
: (J_{Z,W}^{12}(\la^1))^{-1}\,
\sum_{U \in Ir} (\tau^U_{Z,W})^{12}\,(\Id_{H^U_{Z,W}}\T L^U)\,\times
\notag
\\
(\ta _{Z,W}^U)^{12}\,J^{12}_{Z,W}(\la^2):\,=
\sum_{ij}(E^V_{ij})^3:
(J_{Z,W}^{12}(\la^1-\om_i))^{-1}\,(L^V_{ij})^4
\sum_{U \in Ir} (\tau^U_{Z,W})^{12}\,\times
\notag
\\
(\Id_{H^U_{Z,W}}\T L^U)\,
(\ta _{Z,W}^U)^{12}\,J^{12}_{Z,W}(\la^2-\om_j):\,=
\notag
\\
: (J_{Z,W}^{12}(\la^1- h^{(3)}))^{-1}\,
\sum_{U \in Ir} (\tau^U_{Z,W})^{12}\,(L^V)^{34}(\Id_{H^U_{Z,W}}\T L^U)\,
(\ta _{Z,W}^U)^{12}\,J^{12}_{Z,W}(\la^2-h^{(3)}):\,=
\notag
\\
: (J_{Z,W}^{12}(\la^1- h^{(3)}))^{-1}\,
\sum_{U \in Ir} (\tau^U_{Z,W})^{12}\,
(J_{U,V}^{23}(\la^1))^{-1}\,
\notag
\\
\sum_{Y \in Ir} (\tau^Y_{U,V})^{23}\,
(\Id_{H^Y_{U,V}}\T L^Y)\,
(\ta _{U,V}^Y)^{23}\,J^{23}_{U,W}(\la^2)
(\ta _{Z,W}^U)^{12}\,J^{12}_{Z,W}(\la^2-h^{(3)}):\,=
\notag
\\
: (J_{Z,W}^{23}(\la^1- h^{(3)}))^{-1}\,
(J_{Z\T W,V}^{12,3}(\la^1))^{-1}\,
\sum_{U \in Ir} (\tau^U_{Z,W})^{12}\,
\notag
\\
\sum_{Y \in Ir} (\tau^Y_{U,V})^{23}\,
(\Id_{H^Y_{U,V}}\T L^Y)\,
(\ta _{U,V}^Y)^{23}\,
(\ta _{Z,W}^U)^{12}\,
J^{12,3}_{Z\T W, V}(\la^2)\,
J^{12}_{Z,W}(\la^2-h^{(3)}):\,.
\eean
Formulas \Ref{RHS} and \Ref{LHS} and Theorem \ref{R}
imply the Theorem.
\end{proof}

\begin{thm} 
For $V,W \in Ir$, we have
\bean\label{rmr}
R^{12}_{V,W}(\la^1)(L^V)^{13}(L^W)^{23}=:(L^W)^{23}(L^V)^{13}R^{12}_{V,W}(\la^2):
\eean
where the normal ordering sign :: as before means that the matrix elements of $L$
should be put on the right of the matrix elements of $R$. Thus, if
 $L^V=\sum E_{ij}\T L^V_{ij}$, $L^W=\sum E_{kl}\T L^W_{kl}$,
$R(\la)=\sum E^V_{ij}\T E^W_{kl}\T R_{ijkl}(\la)$, then \Ref{rmr} has the form
\bean
\sum_{j,l} R_{ijkl}(\la^1)L^V_{jj'}L^W_{ll'}=
\sum_{j,l} R_{jj'll'}(\la^2)L^W_{kl}L^V_{ij}.
\notag
\eean

\end{thm}
\begin{proof}
\bean
\text{RHS}=:(L^W)^{23}(L^V)^{13}R^{12}_{V,W}(\la^2):=
\notag
\\
: (J_{V,W}^{12}(\la^1))^{-1}\,
\sum_{U \in Ir} (\tau^U_{V,W})^{12}\,(\Id_{H^U_{V,W}}\T L^U)\,
(\ta _{V,W}^U)^{12}\,J^{12}_{V,W}(\la^2) R^{12}_{V,W}(\la^2):=
\notag
\\
: (J_{V,W}^{12}(\la^1))^{-1}\,
\sum_{U \in Ir} (\tau^U_{V,W})^{12}\,(\Id_{H^U_{V,W}}\T L^U)\,
(\ta _{V,W}^U)^{12}\, \Rm^{21}|_{V\T W}P_{W,V}P_{V,W} J^{21}_{W,V}(\la^2):
\notag
\eean
Since $\Rm^{21}|_{V\T W}P_{W,V}$ is an intertwiner, the last expression is equal to
\bean
: (J_{V,W}^{12}(\la^1))^{-1}\,\Rm^{21}|_{V\T W} P_{W,V}
\sum_{U \in Ir} (\tau^U_{W,V})^{12}\,(\Id_{H^U_{W,V}}\T L^U)\,
(\ta _{W,V}^U)^{12}\, P_{V,W} J^{21}_{W,V}(\la^2):=
\notag
\\
: (J_{V,W}^{12}(\la^1))^{-1}\,\Rm^{21}|_{V\T W} P_{W,V}
J^{12}_{W,V}(\la^1)(L^V)^{23}(L^W)^{13} P_{V,W}:=
\notag
\\
: (J_{V,W}^{12}(\la^1))^{-1}\,\Rm^{21}|_{V\T W} 
J^{21}_{W,V}(\la^1)(L^V)^{13}(L^W)^{23}:= \text{LHS}\,.
\notag
\eean
\end{proof}

We proved that $E$ is an $A_0$-algebra. Hence
$E\Tb E$ is an $A_0$-algebra. Define {\it
a comultiplication} $\Dl : E \to E \Tb E$ by
\bean
\Dl f(\la^1)= f(\la^1)\,
\qquad
\Dl f(\la^2)= f(\la^2)\,
\qquad
\Dl (L^V)\,=\,(L^V)^{12}\,(L^V)^{13}
\notag
\eean
where $\Dl (L^V)$ means that $\Dl$ acts in the second factor.
\begin{thm}
The map $\Dl$ preserves the defining relations in $E$.
\end{thm}
\begin{proof} Relations \Ref{mmm} are obviously preserved. 
We check that relation \Ref{dr} is preserved. Compute the image
under $1\T 1\T \Dl$ of the LHS and RHS of \Ref{dr}.
The elements $(1\T 1\T \Dl)\,\text{LHS}$, $(1\T 1\T \Dl)\,\text{LHS}$,
lie in $W\T V \T E \Tb E$. Denote $\la^1_1,\la^2_1$ the $\la$-variables
of the third factor, and
$\la^1_2,\la^2_2$ the $\la$-variables
of the fourth. We have
\bean
(1\T 1\T \Dl)\,\text{LHS}=(L^V)^{23}(L^V)^{24}(L^W)^{13}(L^W)^{14}=
(L^V)^{23}(L^W)^{13}(L^V)^{24}(L^W)^{14}=
\notag
\\
: (J_{W,V}^{12}(\la^1_1))^{-1}\,
\sum_{U \in Ir} (\tau^U_{W,V})^{12}\,(\Id_{H^U_{W,V}}\T L^U)^{123}\,
(\ta _{W,V}^U)^{12}\,J^{12}_{W,V}(\la^2_1):\,\times
\notag
\\
\label {ah}
: (J_{W,V}^{12}(\la^1_2))^{-1}\,
\sum_{Y \in Ir} (\tau^Y_{W,V})^{12}\,(\Id_{H^Y_{W,V}}\T L^Y)^{124}\,
(\ta _{W,V}^Y)^{12}\,J^{12}_{W,V}(\la^2_2):\,.
\eean
We cancel $J^{12}_{W,V}(\la^2_1)$ and $ (J_{W,V}^{12}(\la^1_2))^{-1}$
since in $E\Tb E$ we have a relation $f(\la^2_1)a \Tb b=
a\Tb f(\la^1_2)b$. We replace 
$\ta _{W,V}^U\,\sum_{Y \in Ir} \tau^Y_{W,V}$
with $\Id_{H^U_{W,V}}\T \Id_U$ and use the relation
$f(\la^2)(a \Tb b)=
a\Tb f(\la^2_2)b$ in 
$E\Tb E$.
Thus,
\bean
(1\T 1\T \Dl)\,\text{LHS}=
: (J_{W,V}^{12}(\la^1_1))^{-1}\,
\sum_{U \in Ir} (\tau^U_{W,V})^{12}\,(\Id_{H^U_{W,V}}\T L^U)^{123}\,
\notag
\\
(\Id_{H^U_{W,V}}\T L^U)^{124}\,
(\ta _{W,V}^U)^{12}\,J^{12}_{W,V}(\la^2_2):\,=
(1\T 1\T \Dl)\,\text{RHS}.
\notag
\eean
\end{proof}

For $V\in \Oc_0$,  define $\Id_{V[\mu]}:V\to V$  by 
$\Id_{V[\mu]}|_{V[\mu]}=\Id$ and
$\Id_{V[\mu]}|_{V[\nu]}=0$ for $\nu \neq \mu$.

Define {\it a counit} $\epe : E \to D_T$ where $D_T$ is the $A_0$-algebra
of scalar difference operators on $T$. Set
\bean\label{counit}
\epe (L^V) = \oplus_\mu \Id_{V[\mu]}\T \Ti^{-1}_\mu,
\qquad \epe(f(\la^j))=f(\la).
\eean
\begin{thm}
$\epe$ is a counit in $E$.
\end{thm}
\begin{proof}
The relation 
$$
(\epe \T 1)\,\Dl= (1\T \epe )\Dl = \Id
$$ 
is obviously true.

We  check that the counit $\epe$ preserves the relation \Ref{dr}.
We have
\bean
\epe(LHS)=\epe ((L^V)^{23}(L^W)^{13})\,=\,
\oplus_{\mu,\nu}\Id_{W[\mu]}\T \Id_{V[\nu]}\T \Ti^{-1}_{\mu+\nu}.
\notag
\eean
\bean\label{r-epe}
\epe (RHS)\,=\,
: (J_{W,V}^{12}(\la^1))^{-1}\,
\sum_{U \in Ir} (\tau^U_{W,V})^{12}\,
\sum_{\theta} (\Id_{H^U_{W,V}}\T \Id_{U[\theta]}\T \Ti^{-1}_\theta)\,
(\ta _{W,V}^U)^{12}\,J^{12}_{W,V}(\la^2):\,.
\eean
Notice  that 
$$
\oplus_{U\in Ir} \tau^U_{W,V}
(\Id_{H^U_{W,V}}\T \Id_{U[\theta]})\ta _{W,V}^U=
\Id_{(W\T V)[\theta]}.
$$
Returning to \Ref{r-epe} we get
\bean
\epe (RHS)\,=\,
: (J_{W,V}^{12}(\la^1))^{-1}\,
\sum_{\theta} (\Id_{(W\T V)[\theta]}\T \Ti^{-1}_\theta)
\,J^{12}_{W,V}(\la^2):\,=
\notag
\\
\sum_{\theta} \Id_{(W\T V)[\theta]}\T \Ti^{-1}_\theta =\epe (LHS).
\notag
\eean
The Theorem is proved.
\end{proof}

\subsection {The antipode in $E$}\label{ANT}

\begin{lemma}
If $S : E\to E$ is an antipode, then $S(L^V)=(L^V)^{-1}$,
where $(L^V)^{-1}\in \End (V)\T E$ is such that

$L^V (L^V)^{-1}=\Id_V\T 1$ and $(L^V)^{-1}L^V=\Id_V\T 1$.
\end{lemma}
\begin{proof}
The axioms of the antipode are 
\bean
m\circ (\Id \T S) \circ \Dl (x)=\mu_l(\epe(x)\cdot 1),
\qquad
m\circ (S\T \Id) \circ \Dl (x)=\mu_r(\epe(x)\cdot 1).
\notag
\eean
Applying the first axiom to $L^V$ we get
\bean
\text{LHS}\,:\,L^V\to (L^V)^{12}(L^V)^{13}
\buildrel 1\T S \over \longrightarrow
(L^V)^{12}\,S(L^V)^{13}
\buildrel m \over \longrightarrow
L^V \,S(L^V),
\notag
\\
\text{RHS}\,:\,
L^V\to \sum_\theta \Id_{V[\theta]}\T \Ti^{-1}_\theta \to
\Id_V\T 1.
\notag
\eean
Thus, $L^V S(L^V)=\Id_V\T 1$. Similarly, applying the second axiom, 
we get $S(L^V)L^V=\Id_V\T 1$.
\end{proof}

For $V\in \Oc_0$ define operators $\tK (\la) : {}^*V\to {}^*V$ and 
$K'(\la) : {}^*V\to {}^*V$ by 
\bean\label{K}
\tK (\la)=m(J^{t_2}_{{}^*V,V}(\la)), \qquad
K' (\la)=m(J^{t_1}_{V,{}^*V}(\la)),
\eean
where ${}^{t_j}$ means the  dualization in the $j$-th component,
$(\sum a_i\T b_i)^{t_1}=\sum a_i^*\T b_i$, and $m(a\T b)=ab$. 

 If $\tK(\la)$ is invertible, then denote $K(\la)=(\tK(\la - h))^{-1}$.
Set
\bean\label{lb}
\bar{L}^V= (: K^{(1)}(\la^1)L^{{}^*V}(K^{(1)}(\la^2))^{-1}:)^{t_1},
\\
\label{lh}
\hat{L}^V= (: {K'}^{(1)}(\la^1)L^{{}^*V}({K'}^{(1)}(\la^2))^{-1}:)^{t_1}.
\eean

\begin{thm}\label{antip}
Suppose that $\tK$ or $K'$ is invertible for any module $V\in Ir$.
Then $E=E(\Oc_0)$ is an $A_0$-Hopf algebroid with the antipode
$S(f(\la^1))=f(\la^2),\,S(f(\la^2))=f(\la^1)$ and
$S(L^V)=(L^V)^{-1}= \Lb^V = \hat {L}^V$. Moreover, $K=K'$.
\end{thm}

The Theorem is proved by direct verifications.

The $A_0$-Hopf algebroid $E(\Oc_0)$ will be called {\it the exchange dynamical
quantum group associated to the category $\Oc_0$}.

\subsection{The two point function and $K'(\la)$}\label{tpt}
Define a bilinear from $B_{\la,V} : V\T {}^*V \to \C$. 
For homogeneous $v\in V, v^*\in  {}^*V$, with weights $\la_{v}+\la_{v^*}\neq 0$
set $B_{\la,V} (v,v^*)=0$. If $\la_{v}+\la_{v^*}=0$, then define
$B_{\la,V} (v,v^*)$ by the property
\bean\label{B}
(1\T <\,,\,>_{V\T{}^*V})\circ \Phi^{v,v^*}_\la=
B_{\la,V} (v,v^*)\, \Id_{M_\la}.
\notag
\eean
Notice that $(1\T <\,,\,>_{V\T{}^*V})\circ \Phi^{v,v^*}_\la$ is an intertwiner,
hence it has the form: Const $ \Id_{M_\la}$.
The bilinear form $B_{\la,V}$ is called {\it the two point function}.

\begin{lemma}
$B_{\la,V} (v,v^*)=<v, K'(\la) v^*>$ where $K'(\la)$ is defined
in \Ref{K}.

%=\sum a_i^*b_i$ if $J_{V,{}^*V}(\la)=\sum a_i\T b_i$.

\end{lemma}\label{2pt}
\begin{proof}
Since $\Phi^{v,v^*}_\la=\Phi^{J_{V,{}^*V}(\la)(v\T v^*)}_\la$,
we have
$B_{\la,V} (v,v^*)=\sum < a_iv, b_i v^*>= <v, K'(\la)v^*>$.
\end{proof}

{\bf Remark.} Let $k,n$ be natural numbers, $U$  the vector
representation of the quantum group $U_q(sl_n)$. Let
$V=S^{kn}U$  be the $kn$-th symmetric power of $U$. Then
$V[0]$ is one dimensional and $B_{\la,V}|_{V[0]}$ is a scalar function
of $la$ equal to the squared norm of
a Macdonald polynomial, see Theorem 2.4 in [EK].

\section{Exchange quantum groups associated to simple Lie algebras}

\subsection{The exchange quantum groups $F(\g), F_q(\g)$}
\label{5.1}
In this Section  we consider the exchange dynamical quantum 
groups associated to the category of  finite dimensional representations
of simple Lie algebras and their quantum groups.
We consider two types of polarized Hopf algebras.
\begin{enumerate}
\item[I.] Let $\g$ be a simple Lie algebra, $\al_i, i=1,...,r,$ simple roots, 
$e_i, f_i, h_i$ the corresponding Chevalley generators,
$ \g=\n_+\oplus\h\oplus \n_- $
the polar decomposition. Consider the polarized Hopf algebra $A=U(\g)$
with the $\Z$-grading and polarizations defined by
deg $(e_i)=1$, deg $(f_i)=-1$, deg $(h_i)=0$,
$A_+=U(\n_+)$, 
$A_-=U(\n_-),{}\, 
A_0=U(\h), {} \,$ 
$ A_{\geq 0}=U(\frak b_+)$,
$A_{\leq 0}=U(\frak b_-)$, where $\frak b_{\pm}=\h\oplus \n_{\pm}$. In this case
$T=$ Spec $A_0=\h^*$. Fix on $A$ the quasitriangular structure 
 $\Rm = 1 \in A\,\hat {\T}\, A$.
\item[II.]
Fix $\eps \in\C$ and set $q=e^\eps$. Assume that $q$ is not a root of unity.
Let $\g$ be a simple Lie algebra, $\al_i, i=1,...,r,$ simple roots, 
$\g=\n_+\oplus\h\oplus \n_-$ the polar decomposition. 
Consider the quantum group $A=U_q(\g)$ with the Chevalley generators
$e_i, f_i, K_i^{\pm 1}$ as defined on p. 280 in [CP]. Fix in $A$
a counit $\epe$, a comultiplication $\Dl$, and an antipode
$S$ as defined on p. 281 in [CP]. We consider $A$ as a polarized
Hopf algebra with the $\Z$-grading and polarizations defined by
deg $(e_i)=1$, deg $(f_i)=-1$, deg $(K^{\pm 1}_i)=0$,
$A_+=U_q(\n_+)$, $A_-=U_q(\n_-)$, $A_0=U_q(\h)$, $A_{\geq 0}=U_q(\frak b_+)$,
$A_{\leq 0}=U_q(\frak b_-)$.

{\bf Remark.}
Let $a_{ij}=2<\al_i,\al_j>/<\al_i,\al_i>$ be the Cartan matrix. Let $d_i$ be
coprime positive integers such that the matrix $d_ia_{ij}$ is symmetric.
Let $h_i \in \h$ be the elements
such that $\al_i(h_j)=a_{ij}$. Then one can think of
the generators $K_i^{\pm 1}$ as of elements of the form $q^{\pm d_ih_i}$,
see p. 281 in [CP].

For $A=U_q(\g)$, the spectrum $T=$ Spec $A_0$ is the spectrum of the algebra
of Laurent polynomials $\C[K_1^{\pm 1}, ..., K_r^{\pm 1}]$. The spectrum $T$
can be identified with $\h^*/L$ where $L$ is the lattice such that its dual lattice
$L^*$ is generated by elements $d_ih_i$, i.e. the lattice $L$ consists of the points
where $q^{d_ih_i}$ are equal to $1$. 
%The lattice $L$ will be denoted $P^\vee$.

Fix on $A$ the quasitriangular structure $\Rm  \in A\,\hat {\T}\, A$
where $\Rm$ is the universal R-matrix of the quantum group $U_q(\g)$.
\end{enumerate}

{\bf Remark.} If $q=1$, then sometimes we shall use the notation $U_{q=1}(\g)$ for
the universal enveloping algebra $U(\g)$ considered above.

If $A=U(\g)$, then let $\Oc_0$ be the category of finite dimensional modules over $U(\g)$.
If $q\neq 1$ and $A=U_q(\g)$, then let $\Oc_0(q)$ be the category of 
finite dimensional modules over 
$U_q(\g)$ such that all of eigenvalues of $K_i$ are
integer powers of $q$, i.e. $\Oc_0(q)$ is the category of
finite dimensional modules over $U_q(\g)$ which are quantizations of finite dimensional
modules of $U(\g)$ when $q$ tends to $1$.

Consider the exchange dynamical quantum group $E(\Oc_0)$ associated to
 the category $\Oc_0$ of modules over $U(\g)$ and denote it $\hat F(\g)$.
The  exchange dynamical quantum group $\hat F(\g)$ is defined over
the field of meromorphic functions $M_T,\, T=\h^*$. We replace the field of meromorphic
functions $M_T$ by the field of rational functions $\C(T)$
and define in the same way
{\it the rational   exchange dynamical quantum group} $E_{\text{rat}}(\Oc_0)$.
We denote the rational exchange dynamical quantum group $F(\g)$.

If $q\neq 1$, then
consider   the exchange dynamical quantum group $E(\Oc_0(q))$ associated to
the category $\Oc_0(q)$ of modules over $U_q(\g)$ and denote it $\hat F_q(\g)$.
The exchange dynamical quantum group $\hat F(\g)$ is defined over
the field of meromorphic functions $M_T$, where the torus $T$ has the form 
$T=\h^*/L$. We replace the field of meromorphic
functions $M_T$ by the field of rational functions $\C(T)$
and define in the same way {\it the rational
exchange dynamical quantum group} $E_{\text{rat}}(\Oc_0(q))$.
The field $\C(T)$ can be considered as the subfield $\C(T)\subset M_{\h^*}$
of "trigonometric" functions with respect to the lattice $L\subset
\h^*$. We denote the rational exchange dynamical quantum group $F_q(\g)$.

According to Theorem \ref{antip}, the exchange quantum group $F(\g)$ (resp.
$F_q(\g)$) has a well defined 
antipode if for any $V\in Ir\subset \Oc_0$ (resp. $V\in Ir\subset
\Oc_0(q)$ )
the operator $K'(\la):{}^*V\to {}^*V$ is invertible
for generic values of $\la$. By Lemma \ref{2pt} this property holds if the two point
function $B_{\la,V}$ is a nondegenerate bilinear form for generic values of $\la$.

\begin{thm}\label{nond2pt}
For any $V\in \Oc_0$ (resp. $V\in \Oc_0(q)$ for generic $q$)
the two point function $B_{\la,V} :V\T {}^*V\to \C$ 
is a nondegenerate bilinear form for generic values of $\la$.
\end{thm}

\begin{proof}
For $F(\g)$ the Theorem follows from the next Lemma.

Recall that $B_{\la,V} (v,v^*)=\sum < a_iv, b_i v^*>$
if $J_{V,{}^*V}(\la)=\sum a_i\T b_i$. Let $\rho \in \h^*$ be
the half sum of positive roots.

\begin{lemma}\label{} 
 For $A=U(\g)$ and any $V,W\in \Oc_0$, we have
$J_{V,{}W}(t\rho)\to 1$ when $t\in \C$ and $t$ tends to infinity.
\end{lemma}
\begin{proof} 
In [ES1], the intertwining operator 
$\Phi^v(\la)$ was computed in terms of the Shapovalov form 
(formula (3-5) in [ES1]). From formula (3-5) in [ES1] it is easy to obtain 
the following asymptotic expansion of $\Phi^v(\la)$:  
\bean
\Phi^w_\la v_\la=v_{\la-wt(w)}\T w +
O({1\over |\la|}),
\notag
\eean
where $O({1\over |\la|})$ denotes terms of degree -1 and lower 
in  $\la$. This implies the Lemma.

\end{proof}

\begin{corollary}
$B_{t\rho,V}(\,,\,)\to <\,,\,>$ as $t$ tends to infinity along the imaginary axis.
\end{corollary}

For $F_q(\g)$ and $|q|<1$ or $|q|>1$ the Theorem follows in a
similar way from [ES2], Chapter 2. However, in the q-case,
the above Lemma holds only for $t\to +\infty$ if $|q|<1$ and 
for $t\to -\infty$ if $q>1$. 

\end{proof}

\subsection {The exchange groups and $A_0$-bialgebroids associated with
R-matrices } 
Let $V\in Ir \subset \Oc_0$ (resp. $V\in Ir\subset \Oc_0(q)$ ). Let
$R(\la)=R_{V,V}(\la): V\T V \to V\T V$ be the R-matrix defined in \Ref{r}.
$R$ is a rational function of $\la\in T$. Consider the rational
$A_0$-bialgebroid $A_{rat,R}$ constructed in Section \ref{hh}.
Recall that $A_{rat,R}$ is generated by matrix elements of operators
$L^{\pm 1}$ and rational functions of $\la^1, \la^2 \in T$.

\begin{thm}\label{homomor} For any $V\in Ir\subset \Oc_0$ (resp. $V\in Ir\subset \Oc_0(q)$ ),
there exists a unique homomorphism 
$\phi:A_{rat,R} \to F (\g)$ (resp.
$\phi:A_{rat,R} \to F_q (\g)$)
of rational $A_0$-bialgebroids
such that $(1\T \phi)(L)=L^V$. Moreover, $(1\T \phi)(L^{-1})=(L^V)^{-1}$,
$\phi (f(\la^1))=f(\la^1)$, $\phi (f(\la^2))=f(\la^2)$.
\end{thm}
The Theorem follows from definitions. $\square$

\begin{thm} For $V\in Ir \subset \Oc_0$ (resp. $V\in Ir \subset \Oc_0(q)$),
let $V$ and ${}^*V$ generate the tensor category $\Oc_0$
(resp. $\Oc_0(q)$) in the sense that any object in $Ir$
 is a sub-object in $V^{\T n}\T ({}^*V)^{\T m}$ for suitable
$n, m$. Then the homomorphism $\phi$ is surjective.
\end{thm}
\begin{proof}
Clearly the matrix components of $L^V$ and $L^{{}^*V}$ belong to the image
of $\phi$, since $(L^V)^{-1}$ is $L^{{}^*V}$ up to some invertible factors 
in $\la^1, \la^2$.

Let $U\in  Ir$ and $U$ is a sub-object in $V^{\T n}\T ({}^*V)^{\T m}$ for suitable
$n, m$. Consider the product
$$
(L^{{}^*V})^{m+n,m+n+1}\cdot ... \cdot (L^{{}^*V})^{n+1,m+n+1}\cdot
(L^{V})^{n,m+n+1}\cdot ... \cdot (L^{V})^{1,m+n+1}.
$$
It is clear that the matrix components of $L^U$ are linear
combinations of the matrix components of this product with
coefficients in rational functions of $\la^1, \la^2$.
\end{proof}

\subsection{The exchange groups corresponding to classical Lie groups
$GL(N),$ $ SL(N),$ $ O(N), $ $SP(2N)$}\label{class}
In this Section we  modify the construction of Section \ref{5.1}.

Consider the Lie algebra $gl(N)$. Let $e_i,f_i,\,i=1,...,N-1,$
and $h_i,\,i=1,...,N,$ be its standard Chevalley generators.
Let $\Oc_0(GL(N))$ be the category of all finite dimensional modules over $gl(N)$
which can be integrated to a representation of the Lie group $GL(N)$.
Consider the rational exchange dynamical quantum group $E_{rat}(\Oc_0(GL(N)))$ associated to
 the category $\Oc_0(GL(N))$  and denote it $F(GL(N))$.

Fix $\eps \in\C$ and set $q=e^\eps$. Assume that $q$ is not a root of unity.
Consider the quantum group $A=U_q(gl(N))$ with the standard
Chevalley generators
$e_i, f_i,\,i=1,...,N-1,$ and $k_i^{\pm 1},\, i=1, ...,N$.
%Let $K_i=k_{i+1}/k_i,\, i=1,...,N-1$.
Let $\Oc_0(GL(N),q)$ be the category of all finite dimensional modules over $U_q(gl(N))$
which are q-deformations of finite dimensional modules over $GL(N)$.
Consider the rational exchange dynamical quantum group $E_{rat}(\Oc_0(GL(N),q))$ associated to
the category $\Oc_0(GL(N),q)$  and denote it $ F_q(GL(N))$.

Similarly, let $G$ be a simple complex Lie group and $\g$ its Lie algebra.
Consider the category $\Oc_0(G)$ of all finite dimensional modules over $\g$
which can be integrated to a module over $G$.
Consider the rational exchange dynamical quantum group 
$E_{rat}(\Oc_0(G))$ associated to
 the category $\Oc_0(G)$  and denote it $ F(G)$.
If $\eps \in\C$, $q=e^\eps$, and  $q$ is not a root of unity,
consider the quantum group $A=U_q(\g)$ and the category
$\Oc_0(G,q)$ of all finite dimensional modules over $U_q(\g)$
which are q-deformations of finite dimensional modules over $G$.
The rational exchange dynamical quantum group $E_{rat}(\Oc_0(G,q))$ associated to
the category $\Oc_0(G,q)$ is denoted  $ F_q(G)$.

Let $G$ be a Lie group of type $GL(N), SL(N), SO(N), SP(2N)$
and $\g$ its Lie algebra.
Let $V$ be the vector representation of $U(\g)$ (resp. $U_q(\g)$).
We have $V\in Ir\subset \Oc_0(G)$ (resp. $V\in Ir\subset \Oc_0(G,q)$).

\begin{lemma} 
$V$ and ${}^*V$ generate $\Oc_0(G)$ (resp. $\Oc_0(G,q)$).
\end{lemma}
The Lemma follows from the fact that the vector representation is faithful
as a representation of $G$.

\begin{corollary}\label{vect-rep}
Let $V$ be the vector representation of $U(\g)$ (resp. $U_q(\g)$),
$R(\la)=R_{V,V}(\la): V\T V \to V\T V$ the R-matrix defined in \Ref{r},
$A_{rat,R}$ the rational $A_0$-bialgebroid  constructed in Section \ref{hh}.
Then the homomorphism $\phi:A_{rat,R} \to F (G)$ (resp.
$\phi:A_{rat,R} \to F_q (G)$) of Theorem \ref{homomor} is an epimorphism.
\end{corollary}

\begin{thm} \label{gln}
Let $G=GL(N)$. Then
\begin{enumerate}
\item[I.] For $F (G)$,   the homomorphism $\phi:A_{rat,R} \to F (G)$  
of Corollary \ref{vect-rep} is injective.
\item[II.] For $F_q (G)$,   the homomorphism $\phi:A_{rat,R} \to F_q (G)$  
of Corollary \ref{vect-rep} is injective for all $q$ except a countable
set.
\end{enumerate}
\end{thm}
\begin{proof}
To prove the Theorem for $F (G)$ recall  that in this case $\la\in T=\h^*$.
For $\gm\in \C^*$ introduce a new variable $\tilde \la=\la/\gamma$. Then, 
by the results of Chapter 3 in [ES1], for any modules
$V,W \in \Oc(GL(N))$, we have
$J_{V,W}(\lt )= \Id + \gm J_1(\lt ) + \gm^2 J_2 (\lt ) + ... $.
Hence $J_{V,W}(\lt )\to \Id$ as $\gm \to 0$.

Let $A_{rat,R}^\gm, F^{\gm} (G)$   be the algebras defined by
the same relations as $A_{rat,R}, F (G)$
with $\la$ replaced by $\la/\gm$. It is easy to see that
the algebras $A_{rat,R}^0, F^{0} (G)$ 
are well defined,
$$
A_{rat,R}^0 = F^{0} (G)= \C(\h^*)\T \C(\h^*) \T \C[G]
$$
and $\phi_\gm \to \phi_0=\Id $ as $\gm\to 0$.
Here $\C[G]$ is the algebra of polynomials on $G$.

The algebras $A_{rat,R}^\gm, F^{\gm} (G)$ and the homomorphism
$\phi_\gm$  are deformations
of the algebras $A_{rat,R}^0, F^{0} (G)$ and the homomorphism
$\phi_0$. Elementary reasonings of the deformation theory imply that
the homomorphism $\phi_\gm$ is an isomorphism.

The Theorem for $F_q (G)$ is deduced from the Theorem for $F (G)$
by taking the limit $q\to 1$.
\end{proof}

Now consider the case of $SL(N)$. 

For $G=GL(N)$, 
consider the exchange group $F_q(G)$. Let $C\in \Oc_0(G,q)$
be a one dimensional module. Then $L^C$ is a $1\times 1$-matrix and
can be considered as an element of $F_q(G)$.

\begin{lemma}
$L^C$ is a central element in $F_q(G)$ and $L^C$ is invertible,
$(L^C)^{-1}=L^{{}^*C}$.
\end{lemma}
\begin{proof} For any $W\in \Oc_0(G,q)$, $L^C$ and $L^W$ satisfy
the R-matrix relation \Ref{rmr}. In this case the R-matrix $R_{C,W}(\la)$
is a scalar constant, hence $L^C$ is central.
\end{proof}

For $F_q(GL(N))$, consider the one dimensional module
$C=\wedge^N_q V$ over $U_q(gl (N)),$
 which is the $N$-th quantum exterior power of the vector representation
$V$. For generic $q$, consider the isomorphism $\phi_{GL(N)}
: A^{GL(N)}_{rat,R} \to F_q (GL(N))$  of Theorem \ref{gln}.
Define $D \in A^{GL(N)}_{rat,R}$ by $D=\phi_{GL(N)}^{-1}(L^C)$.

Consider the quantum group $U_q(sl(N))$. There is a natural embedding
of $U_q(sl(N))$ to $ U_q(gl(N))$ sending the Chevalley generators $e_i,f_i, K_i
\in U_q(sl(N))$ to $e_i,f_i, k_{i+1}/k_i \in U_q(gl(N))$.
Let $V$ be the vector representation of $ U_q(sl(N))\subset
U_q(gl(N))$. Consider the corresponding R-matrices $R^{GL(N)}(\la)=
R^{GL(N)}_{V,V}(\la), \, \la \in T_{GL(N)}=(\C^*)^N$ and $R^{SL(N)}(\la)=
R^{SL(N)}_{V,V}(\la), \, \la \in T_{SL(N)}=(\C^*)^N/\C^*(1,...,1)$.
Any rational function on $(\C^*)^N/\C^*(1,...,1)$ can be considered as a rational
function on $(\C^*)^N$ invariant with respect to the diagonal action of $\C^*$.
It is easy to see that R-matrix
$R^{SL(N)}(\la)$ considered as a function on $T_{GL(N)}$
coincides with the R-matrix $R^{GL(N)}(\la)$ up to a multiplicative
scalar constant. This construction allows us to
define a natural embedding $A^{SL(N)}_{rat,R} \to A^{GL(N)}_{rat,R}$.
Clearly, the element $D$ belongs to the image of the imbedding.

\begin{thm}

\begin{enumerate}
\item[I.] For $F (SL(N))$,   the kernel of the epimorphism $\phi:A^{SL(N)}_{rat,R} 
\to F (SL(N))$ of Corollary \ref{vect-rep} is generated by the relation
$D=1$.
\item[II.] For $F_q (SL(N))$,  the kernel of the epimorphism $\phi:A^{SL(N)}_{rat,R}
 \to F_q (SL(N))$  
of Corollary \ref{vect-rep} contains the ideal  generated by  the relation $D=1$.
Moreover, the kernel is generated by  this relation
for all $q$ except a countable set.
\end{enumerate}
\end{thm}
\begin{proof} For $F (SL(N))$,  clearly the kernel contains the relation $D=1$,
since for $sl(N)$, the module $\wedge^N V$ is trivial. 

Introduce (as before)
the algebras $A^{SL(N),\gm}_{rat,R}, F^{\gm} (SL(N))$  
and a homomorphism $\phi_\gm : A^{SL(N),\gm}_{rat,R} \to F^{\gm} (SL(N))$  depending on
a parameter $\gm \in \C^*$.
It is easy to see that for $\gm= 0$, the homomorphism
$\bar \phi_{\gm=0}:
A^{SL(N),\gm=0}_{rat,R}/\{D=1\} \to F^{\gm=0} (SL(N))$  
is an isomorphism. This statement (as before) implies the Theorem.
\end{proof}

Now let $G$ be a Lie group of type $SO(N)$ or $SP(2N)$. Let $V$ be
the vector representation of its Lie algebra $\g$ (resp. $U_q(\g)$).
In this case there is an isomorphism $T:{}^* V \to V$ of $\g$-modules
(resp. $U_q(\g)$-modules).

\begin{thm} Let $G$ be a Lie group of type $SO(N)$ or $SP(2N)$. Then
\begin{enumerate}
\item[I.] For $F (G)$
and $F_q(G)$ 
 the kernel of the epimorphism $\phi:A_{rat,R} \to F (G)$  
of Corollary \ref{vect-rep} contains the ideal generated by the relations
\bean\label{sp}
L=: T^{(1)}(K^{(1)}(\la^1))^{-1} (L^{-1})^{t_1}K^{(1)}(\la^2)
(T^{(1)})^{-1}:\,
\eean
where $K$ is defined in Section \ref{ANT}.

\item[II.]
For $F(G)$ and $F_q(G)$,
the element $D$ defined above equals 1 modulo \Ref{sp} for $G=SP(2N)$, and 
is a central grouplike element of order 2 modulo \Ref{sp} for $G=SO(N)$. 
 
\item[III.]  For $F(G)$ and $F_q(G)$ with $q$ outside of a countable set, 
the kernel of $\phi$ is generated by  relations \Ref{sp}
in the case of $G=SP(2N)$, and by \Ref{sp} and $D=1$ for $G=SO(N)$.  
\end{enumerate}

\end{thm}
\begin{proof}
\begin{lemma}
Relations \Ref{sp} belong to the kernel. 
\end{lemma}
\begin{proof}
In fact, by Theorem \ref{antip} we have
$$
{L}^{{}^*V}= : (K^{(1)}(\la^1))^{-1} ((L^{V})^{-1})^{t_1}K^{(1)}(\la^2):\,.
$$
Since $T:{}^*V\to V$ is an isomorphism, we have
$$
(T\T 1) L^{{}^*V}(T^{-1}\T 1)= L^V.
$$
\end{proof}

Let $I\subset A_{rat, R}$ be the ideal generated by relations \Ref{sp}.
Consider the quotient $A_{rat, R}/I$ and the homomorphism
$\bar \phi : A_{rat, R}/I \to F_q (G)$. One can prove as for
$GL(N)$ that the homomorphism $\bar \phi$ is an isomorphism for $q=1$
and for generic $q$ if $G=SP(2N)$, and has kernel generated by $D=1$ 
if $G=SO(N)$. 

{\bf Remark.} If $G=SO(N)$, then it is natural to denote 
the quotient $A_{rat, R}/I$ by $F_q(O(N))$

{\bf Remark.} If $q=1$, then in the limit $\gamma\to 0$ we have 
$J=1$. In this case relations \Ref{sp}
take the form
$$
L= (T\T 1) (L^{-1})^{t_1}(T^{-1}\T 1),
$$
which is the defining relation for the orthogonal and symplectic groups.

\end{proof}

\section{The R-matrix $R_{V,V}(\la)$ for the vector representation of $U_q(gl(N))$}

\subsection{Matrices $J_{V,V}(\la)$ and $R_{V,V}(\la)$}
Let $V=\C^N$ be the vector representation of
$A=U_q(gl(N))$. Let $v_j=(0,...,0,1_j,...,0)$ be the standard basis in $V$.
We have $f_iv_j=\dl_{i,j}v_{i+1},\,
e_iv_j=\dl_{i+1,j}v_{i}$ where $f_i, e_i$ are the Chevalley generators of
$U_q(gl(N))$.  Introduce a basis  $E_{ij}$ in $ \End (V)$ by $E_{ij}v_k = \dl_{jk} v_i$.

The $U_q(gl(N))$-module $V\otimes V$ has the weight decomposition,
\bean\label{wght}
V\otimes V \,=\, \oplus_{a=1}^N V_{aa} \, \oplus
\oplus_{a < b} V_{ab}\,,
\eean
where $V_{aa}\, =\, \C\,v_a\otimes v_a$ and $V_{ab}\, =\, \C\,v_a\otimes v_b \oplus
\C\, v_b \otimes v_a $ .

The action of the quasi-triangular structure $\Rm  \in A\hat {\T} A$ on $V\T V$
takes the form
$$
\Rm= q\sum_{a=1}^N E_{aa}\otimes E_{aa}\,+
\,\sum_{a\neq b } \,E_{aa}\otimes E_{bb}\,+
\,\sum_{a< b } \,(q-q^{-1})\,E_{ab}\otimes E_{ba}.
$$

Consider the maps 
$$J(\la)=J_{V,V}(\la): V\T V \to V\T V,
\qquad
R(\la)=R_{V,V}(\la): V\T V \to V\T V
$$
defined in \Ref{dj} and \Ref{r}. 
Here $\la \in T=(\C^*)^N$, if $q\neq 1$, and
$\la \in T=\C^N$, if $q= 1$.  We shall use the  coordinates 
$\la=(q^{\la_1},...,q^{\la_N})$ on $(\C^*)^N$ and the coordinates
$\la=(\la_1,...,\la_N)$ on $\C^N$.

Recall that $R(\la)=J^{-1}\Rm^{21}J^{21}$.

\begin{thm}\label{}
\begin{enumerate}
\item[I.] For $F_q (GL(N))$, we have
\bean
J(\la)= \sum_{a,b}\,E_{aa}\otimes E_{bb}\,+\,
\,\sum_{a< b } \,{q^{-1}-q\over
q^{2(\la_a-\la_b +b-a)}-1}
\,E_{ba}\otimes E_{ab},
\notag
\\ \label{Rtr}
R(\la)= q\sum_{a=1}^N\,E_{aa}\otimes E_{aa}\,+\,
\,\sum_{a\neq b} \,{q^{-1}-q \over
q^{2(\la_b-\la_a +a-b)}-1}
\,E_{ba}\otimes E_{ab}\,+\,
\sum_{a< b} \,E_{aa}\otimes E_{bb}\,+\,
\\
\,\sum_{a> b} \,
{(q^{2(\la_b-\la_a +a-b)}-q^{-2})
(q^{2(\la_b-\la_a +a-b)}-q^{2})
\over
(q^{2(\la_b-\la_a +a-b)}-1)^2}\,
E_{bb}\otimes E_{aa}\,.
\notag
\eean
\item[II.] For $F (GL(N))$, we have
\bean
J(\la)= \sum_{a,b}\,E_{aa}\otimes E_{bb}\,+\,
\,\sum_{a< b } \,{1\over
\la_b-\la_a +a-b}
\,E_{ba}\otimes E_{ab},
\notag
\\ \label{Rrat}
R(\la)= \sum_{a=1}^N\,E_{aa}\otimes E_{aa}\,+\,
\,\sum_{a\neq b} \,{1 \over
\la_a-\la_b +b-a}
\,E_{ba}\otimes E_{ab}\,+\,
\sum_{a< b} \,E_{aa}\otimes E_{bb}\,-\,
\\
\,\sum_{a> b} \,
{(\la_b-\la_a +a-b-1)(\la_b-\la_a +a-b+1)
\over
(\la_b-\la_a +a-b)^2}
E_{bb}\otimes E_{aa}\,.
\notag
\eean

\end{enumerate}
\end{thm}
The Theorem is proved by direct calculations.
More precisely, the coefficients of $J$ corresponding 
to simple roots (i.e. with $b=a+1$) are easily calculated
explicitly, after which all 
 other coefficients are found using the classification 
of dynamical quantum R-matrices of Hecke type given in [EK2]. 

Set $R^\vee(\la)= P R(\la)$ where $P:V\T V\to V\T V$ is the 
permutation of factors. 

\begin{lemma}

\begin{enumerate}
\item[I.] The operator
$R^\vee(\la)$ preserves the weight decomposition \Ref{wght}.
\item[II.] For any $a=1,...,N$, 
we have $R^\vee(\la)v_a\T v_a= q\,v_a\T v_a$.
\item[III.] For any $a\neq b$,
the operator $R^\vee(\la)$ restricted to the two dimensional space
$V_{ab}$ has eigenvalues $q$ and $-p$, where $p=q^{-1}$.
\end{enumerate}
\end{lemma} $\square $

A meromorphic function $R: T \to \End(V\T V)$ with these three properties 
and satisfying the dynamical Yang-Baxter equation \Ref {qdyb} is called
{\it an R-matrix of Hecke type with parameters $q$ and $p$}. We classified
such R-matrices in [EV2] up to gauge transformations.

\subsection{Gauge transformations of R-matrices of Hecke type}
Concider the torus $T=(\C^*)^N$ with  coordinates 
$\la=(q^{\la_1},...,q^{\la_N})$. {\it A multiplicative  $k$-form} on $T$
 is a collection,
$$
\phi \,=\,\{ \phi_{a_1, ... , a_k}(q^{\la_1}, ... ,
 q^{\la_N})\}\,,
$$
of meromorphic
functions on $T$, where $a_1,..., a_k$ run through
all $k$ element subsets of $\{ 1, ..., N\}$, such that for any
subset $a_1,..., a_k$  and any $i,\, 1\leq i < k$, we have
$$
\phi_{a_1, ..., a_{i+1}, a_i, ...  , a_k}(q^{\la_1}, ..., q^{\la_N})
\, \phi_{a_1, ... , a_k}(q^{\la_1}, ..., q^{\la_N}) \, = \, 1 \,.
$$
Let $\Omega^k$ be the set of all multiplicative $k$-forms.

If $\phi$
and $ \psi$
are  multiplicative $k$-forms, then
$\{ \phi_{a_1, ... , a_k}(q^{\la_1}, ..., q^{\la_N}) \, \cdot \,
\psi_{a_1, ... , a_k}(q^{\la_1}, ..., q^{\la_N})\}$
and
$\{ \phi_{a_1, ... , a_k}(q^{\la_1}, ..., q^{\la_N})\,/\,
 \psi_{a_1, ... , a_k}(q^{\la_1}, ..., q^{\la_N})\}$
are multiplicative $k$-forms. This gives an abelian group structure on $\Omega^k$.
The zero element in $\Omega^k$ is the form
$\{ \phi_{a_1, ... , a_k}(q^{\la_1}, ..., q^{\la_N}) \equiv 1\}$.

 For any $a = 1, ...., N$, introduce an endomorphism
$\dl_a$
of the 
multiplicative group of nonzero meromorphic functions $f(q^{\la_1}, ..., q^{\la_N})$ on $T$ by
$$
\dl_a \,:\, f(q^{\la_1}, ..., q^{\la_N}) \, \mapsto \,
{f(q^{\la_1}, ..., q^{\la_N})\over f(q^{\la_1}, ..., q^{\la_a}/q, ..., q^{\la_N})}
$$
and a homomorphism $d : \Omega^{k} \to \Omega^{k+1},\,
\phi \mapsto d\phi,$ by
$$
(d \phi)_{a_1, ... , a_{k+1}}(q^{\la_1}, ..., q^{\la_N})\,=\,
\prod_{i=1}^{k+1} \,(\dl_{a_i}\phi_{a_1, ... , a_{i-1},a_{i+1}, ...,a_{k+1}}
(q^{\la_1}, ..., q^{\la_N}))^{(-1)^{i+1}}.
$$
We have $d^2\,=\, 0$ ($0$ means the trivial homomorphism which maps 
everything to the zero element).
A multiplicative form $\phi$ is called {\it closed} if $d\phi = 0$.

Introduce gauge transformations of 
R-matrices, $R \, : \, T  \to \End (V \otimes V)$, of the form 
\bean\label{r-form}
R(\la) \, =\,
\sum_{a = 1}^N \,
\al_{aa}(\la)\,E_{aa}\otimes E_{aa}\,+
\sum_{a \neq b } \,
\al_{ab}(\la)\,E_{aa}\otimes E_{bb}\,+
\,\sum_{a\neq b } \,\bt_{ab}(\la)\,E_{ba}\otimes E_{ab},
\eean
where $\al_{ab}(\la),\bt_{ab}(\la)$ are suitable functions.

\begin{enumerate}

\item[I.]  Let $\{\phi_{ab}\}$ be a meromorphic closed
  multiplicative $2$-form
on $T$.  Set
$$
R(\la) \, \mapsto \,
\sum_{a = 1}^N \,
\al_{aa}(\la)\,E_{aa}\otimes E_{aa}\,+
\sum_{a \neq b } \,\phi_{ab}(\la)\,
\al_{ab}(\la)\,E_{aa}\otimes E_{bb}\,+
\,\sum_{a\neq b } \,\bt_{ab}(\la)\,E_{ba}\otimes E_{ab}.
$$
\item[II.] Let the symmetric group
$S_N$ , the Weyl group of $gl_N$,
act on $T$ and $V$ by permutation of coordinates.
For any permutation $\sigma \in S_N$, set
$$
R(\la) \, \mapsto \, (\sigma\otimes\sigma)\, R(\sigma^{-1} \cdot \lambda)\,
(\sigma^{-1} \otimes \sigma^{-1}) \,.
$$
\item[III.] For a nonzero complex number $c$, set
$$
R( \la ) \, \mapsto \, c\, R(  \la )\,.
$$
\item[IV.] For  an element $\mu \in T$, set
$$
R( \la ) \, \mapsto \, R( \la \, +\,\mu )\,,
$$
(recall that we always use the additive notation for
the standard group structure on $T$).
\end{enumerate}

By Theorem 1.1 in [EV2] any gauge transformation
transforms a matrix satisfying  the QDYB \Ref{qdyb} to
a matrix satisfying  the QDYB \Ref{qdyb}.
In all cases, if the R-matrix is of Hecke type,
then the transformed matrix is of Hecke type.
If the transformation is of type III and the Hecke parameters
of the R-matrix are $q$ and $p$, then the Hecke parameters
of the transformed matrix are $c q$ and $c p$. For all other types
of transformations the Hecke parameters do not change.

Two R-matrices  will be called {\it rationally
equivalent} if one of them can be transformed into another by a sequence
of gauge transformations of types II-IV and of type I with only rational
functions $\phi_{ab}(\la)$. 

The R-matrices of Hecke type were classified in [EV2] up to gauge transformations.
Here are two main examples of that classification.

{\bf Examples.}

\begin{enumerate}
\item[I.] For $q\neq 1$, let the 
R-matrix $R^q : (\C^*)^N\to V\T V$ have the form \Ref{r-form}, where
\bean\label{beta}
\bt_{ab}(\la) \,=\, { q^{-2}\,-\,1 \over 
q^{2( \la_{b}-\la_a)} \,-\,1 }\,,
\notag
\eean
$\al_{aa}=1$ and $\al_{ab}(\la)=\bt_{ab}(\la)+q^{-2}$ for $a\neq b$.
$R^q(\la)$ is an R-matrix of Hecke type with parameters $q$ and $p=q^{-1}$.
\item[II.] For $q=1$, let the R-matrix $R : \C^N \to V\T V$ have the form \Ref{r-form}, where
\bean\label{beta-rat}
\bt_{ab}(\la) \,=\, { 1 \over 
 \la_{a}-\la_b}\,,
\notag
\eean
$\al_{aa}=1$ and $\al_{ab}(\la)=\bt_{ab}(\la)+1$ for $a\neq b$.
$R(\la)$ is an R-matrix of Hecke type with parameters $q=p=1$.
\end{enumerate}

\begin{lemma}

\begin{enumerate}
\item[I.] The R-matrix $R(\la)$  in \Ref{Rtr} is rationally
equivalent to the R-matrix $R^q(\la)$ of the first example.
\item[II.] The R-matrix $R(\la)$  in \Ref{Rrat} is rationally
equivalent to the R-matrix $R(\la)$ of the second example.
\end{enumerate}
\end{lemma}
\begin{proof} To transform $R^q(\la)$ to the R-matrix in \Ref{Rtr} one needs to
make the gauge transformation of type IV with $\mu$ equal to the half sum of
positive roots, then the gauge transformation of type II with $c=q$, and finally,
the gauge transformation of type I corresponding to the closed multiplicative
2-form $\phi_{ab}(\la)$, where
\bean\label{exact}
\phi_{ab}(\la)=
q { q^{2( \la_{a}-\la_b +a-b-1)} \,-\,1 \over 
q^{2( \la_{a}-\la_b +a-b)} \,-\,1 }\, \qquad
\text{for{}}\,{}\, a>b
%q^{-1} + 
%{ q^{-1}\,-\,q \over 
%q^{2( \la_{a}-\la_b +a-b)} \,-\,1 }\, \qquad
%\text{for{}}\,{}\, a>b
\eean
and
$\phi_{ab}(\la)$ for $a<b$ are reconstructed from the multiplicative
"skew symmetry" relation $\phi_{ba}(\la)\phi_{ab}(\la)=1$.
The second statement of the Lemma is proved analogously.
\end{proof}

\subsection{Gauge transformations of $A_0$-bialgebroids corresponding to
R-matrices}
For $q\neq 1$, let $R, \,\tilde R \,:\,  (\C^*)^N\to V\T V$ be two R-matrices 
having the form \Ref{r-form}.
Consider the Lie algebra $gl(N)$ and its Cartan subalgebra $\h$ generated by
elements $h_i$. Consider the polarized algebra $A=U_q(gl(N))$ with the earlier
distinguished
 polarization $A_{\pm },\,A_0=U_q(\h)$.
Let $A_{R}$ and $A_{\tilde R}$ be the 
$A_0=U_q(\h)$-bialgebroids associated to  $R$ and
$\tilde R$, respectively,  and constructed in Section \ref{hh}.
Assume that the R-matrix $\tilde R$ is obtained from the R-matrix
$R$ by a gauge transformation of type II,  III, or IV, then clearly the
$A_0$-bialgebroids $A_{R}$ and $A_{\tilde R}$ are isomorphic.

\begin{thm}
Assume that the R-matrix $\tilde R$ is obtained from the R-matrix 
$R$ by a gauge transformation of type I associated to a multiplicative
2-form $\{\phi_{ab}(\la)\}$ which is  exact, 
$\{\phi_{ab}(\la)\}=d \,\{\xi_a(\la)\}$
where $\{\xi_a(\la)\}$ is a multiplicative 1-form. Then the 
$A_0$-bialgebroids $A_{R}$ and $A_{\tilde R}$ are isomorphic.
\end{thm}
\begin{proof} Let $\xi= \sum \xi_a E_{aa}$. 
If $L$ satisfies the QDYB equation \Ref{ybr}, and $\tilde L$ is such that
$$
L=:\xi^{(1)}(\la^1)\tilde L (\xi^{(1)}(\la^2))^{-1}:\,,
$$
then $\tilde L$ satisfies
$$
\tilde {R}^{12}(\la^1)\tilde {L}^{13}\tilde{L}^{23}=
:\tilde{L}^{23}\tilde {L}^{13}\tilde {R}^{12}(\la^2):
$$
where
$$
\tilde R(\la)= (\xi^{(1)}(\la - h^{(2)}))^{-1} (\xi^{(2)}(\la))^{-1}
R(\la)\xi^{(1)}(\la)\xi^{(2)}(\la- h^{(1)}).
$$
This means that if $R(\la)$ has the form \Ref{r-form}, then $\tilde R(\la)$
is obtained from $R(\la)$ by the gauge transformation of type I corresponding
to the 2-form $\{\phi_{ab}(\la)\}=d \,\{\xi_a(\la)\}$.
\end{proof}

Notice that the multiplicative 2-form given by \Ref{exact} is exact,
$\{\phi_{ab}(\la)\}=d \,\{\xi_a(\la)\}$, where
$$
\xi_a(\la)\,=\,\prod_{b<a}\,q^{-\la_b}\,
(\,q^{2( \la_{a}-\la_b +a-b-1)} \,-\,1\,).
$$
Hence the corresponding bialgebroids are isomorphic.

\section{Elements of representation theory of exchange
groups $F_q(G)$}

\subsection{A construction of representations}
Let $G$ be a simple group and $\g$ its Lie algebra.
Consider an exchange group  $F_q=F_q(G)$.

A dynamical representation $\pi_W:F_q\to D_{T, W}$ is called {\it rational} if
the image of $\pi_W$ consists of difference operators with rational
coefficients.  A homomorphism of dynamical representations
$\phi: W_1\to W_2$ is called {\it rational} if 
the matrix elements of $\phi$ are rational functions.

Denote $\Repq$ the tensor category of rational finite dimensional (dynamical)
representations of $F_q$ and rational morphisms between the representations. 

Let $W\in \Oc_0(G,q)$. Define a rational dynamical representation of $F_q$ on
$W$. Recall that a rational dynamical representation is a diagonalizable
$A_0=U_q(\h)$-module $W$ and a homomorphism of $A_0$-algebras
$\pi_W:F_q \to D_{T,W}$ such that the image of the homomorphism consists
of difference operators with rational coefficients. We consider $W\in \Oc_0(G,q)$ with the
$A_0$-module structure induced by $U_q(\h)\subset U_q(\g)$ and define $\pi_W$ by 
\bean\label{pi}
\pi_W(f(\la^1))= f(\la),\qquad
\pi_W(f(\la^2))=f(\la - h),
\\
\label{pipi}
(1\T \pi^0_W)(L^V)(\la)=R_{V,W}(\la),
\eean
for any $V\in Ir\subset \Oc_0(G,q)$.

Recall that for a dynamical representation $\pi_W:F_q \to D_{T,W}$,
one defines a  map $\pi^0_W: F_q \to \End(W,W\T M_T)$ 
as explained in Section \ref{dfns} and  this map 
uniquely determines $\pi_W$.

\begin{thm}\label{constr}
Formulas \Ref{pi} and \Ref{pipi} define a structure of a rational dynamical
representation of $F_q$ on $W$.
\end{thm}
The Theorem  follows from definitions.

Define a functor $F$ from the category $\Oc_0(G,q)$ of finite dimensional modules
over $U_q(\g)$ (defined in Section \ref{class}) to the category $\Repq$ 
sending an object $W\in \Oc_0(G,q)$ to $(W,\pi_W)$ and sending a morphism
$\al : W\to U$ to the same linear map $\al : W\to U$.

\begin{thm}\label{iso}
\begin{enumerate}
\item[I.]
This construction defines a tensor functor $F : \Oc_0(G,q)\to \Repq$
with a tensor structure
$$
J_{W,U}:F(W)\tilde {\T} F(U)\to F(W\T U),
$$
where $J_{W,U}$ is defined in \Ref{dj}. 
\item[II.] For generic $q$ the map $\Hom_{\Oc_0(G,q)}(W,U)\to
\Hom_{\Repq}(F(W),F(U))$ defined by $F$ is an isomorphism.
Thus, $F$ defines a tensor equivalence of $\Oc_0(G,q)$ onto a full 
subcategory of $\Repq$. 
\end{enumerate}
\end{thm}

{\bf Remark 1.} 
In the next paper we plan to show that $F$ is an equivalence of categories, 
i.e. that any 
object of $\Repq$ is in the image of $F$. 

{\bf Remark 2.}
We see that the representation category of $F_q(G)$ is essentially the 
same as for $U_q(\g)$. A similar result was obtained in \cite{BBB}, 
where it is shown that the quasi-Hopf algebra associated to the dynamical 
R-matrix is twist equivalent to $U_q(\g)$ (for $\g=sl(2)$). 
These two results are closely related, because as follows from \cite{Xu}, 
representation theory of the Hopf algebroid $F_q(G)$ is 
tautologically equivalent to representation theory of the corresponding 
quasi-Hopf algebra. 

{\bf Remark 3.}
Theorem 45 raises a question: why is it interesting to study
dynamical quantum groups if they have the same representation theory as the 
usual ones? In our opinion, it is interesting to study not only tensor 
categories but also their realizations (e.g. tensor functors on them
to other tensor categories), 
which contain extra structure. In particular, $F_q(G)$ and $U_q(G)$ 
correspond to two different realizations of the same 
tensor category. In other words, dynamical quantum groups 
do not provide new tensor categories, but do provide 
new realizations of already known tensor categories. 

Now let us prove the theorem. 
The first part of the Theorem is trivial. The second part of the Theorem
is proved in Section \ref{proof}.
In order to prove the second part
we first prove that any rational morphism $b: F(W)\to F(U)$ does not depend on $\la$
and then  show that there exists $a \in \Hom_{\Oc_0(G,q)}(W,U)$ such that
$b=F(a)$.

\subsection{Rational morphisms}
Let $W\in \Oc_0(G,q)$. Let $\pi_W, \pi^0_W$ be representation maps.
By definition
$$
(1\T \pi^0_W)(L^V)(\la)=R_{V,W}(\la)=
J^{-1}_{V,W}(\la)\Rm^{21}|_{V\T W } J^{21}_{W,V}(\la).
$$
Let $v_i^V$ be a homogeneous basis of $V$, and $v_0^V$ be the
highest weight vector. Let $(v_i^V)^*$ be the dual basis.
Consider the matrix element
$R^{00}_{V,W}(\la)\in \End (W)$ defined by
$$
<y^*,R^{00}_{V,W}(\la)x>=<(v^V_0)^*\T y^*, R_{V,W}(\la)\, v_0^V\T x>,
$$
where $x\in V, y^*\in V^*$.

\begin{lemma}\label{claim}
$R^{00}_{V,W}(\la)$ is a nonzero scalar operator on each weight subspace 
$W[\al], \al\in T$, of $W$. Moroover, the value of the scalar is determined by
$\al$ and does not depend on $W$.
\end{lemma}
\begin{proof}
We have
\bean
<(v^V_0)^*\T y^*, J^{-1}_{V,W}(\la)\Rm^{21}|_{V\T W } J^{21}_{W,V}(\la)\, v_0^V\T x>=
\notag
\\
\label{cl}
<(J^{-1}_{V,W}(\la))^*\,(v^V_0)^*\T y^*, \Rm^{21}|_{V\T W } J^{21}_{W,V}(\la)\, v_0^V\T x>
\eean
Since deg $J=0$ and $J=1 +\sum a_i\T b_i$, deg $a_i<0$, we have
$J^{21}=1 +\sum b_i\T a_i$ and $(J^{-1})^*=1 +\sum c_i\T d_i$,
where deg $c_i<0$. Continuing \Ref{cl}, we get
\bean\label{cl1}
<y^*,R^{00}_{V,W}(\la)x>=
<(v^V_0)^*\T y^*, \Rm^{21}|_{V\T W }\, v_0^V\T x>=
\notag
<(\Rm^{21}|_{V\T W })^*\,(v^V_0)^*\T y^*,  v_0^V\T x>.
\notag
\eean
It is well known that the operator $(\Rm^{21}|_{V\T W })^*$ has the form
$\Rm_0 Q$, where $\Rm_0 =1 +$(a strictly upper triangular element in
$U_q(\n_+)\hat{\T}U_q(\n_-)$)
and $Q\in U_q(\h)\hat{\T} U_q(\h)$. Hence,
$$
<(\Rm^{21}|_{V\T W })^*\,(v^V_0)^*\T y^*,  v_0^V\T x>=
<Q\,(v^V_0)^*\T y^*,  v_0^V\T x>.
$$
This proves the Lemma.
\end{proof}

\begin{lemma}\label{const}
Let $a(\la):F(W)\to F(U)$ be an intertwining operator,
then $a(\la)$ does not depend on $\la$.

\end{lemma}
\begin{proof}
We have $R^{00}_{V,W}=\pi^0_W(L^V_{00})$,
where $L^V_{00}$ is the matrix component of $L^V$ corresponsing
to the highest weight vector $v^V_0$.
Hence, $\pi_W(L^V_{00})=R^{00}_{V,W}\Ti^{-1}_{wt(v_0^V)}$,
where $wt(v_0^V)$ is the weight of $v_0^V$.

The intertwining operator has to satisfy $a(\la)\circ \pi_W(L^V_{00})=
\pi_W(L^V_{00})\circ a(\la)$. Hence, $a(\la)=a(\la - wt(v_0^V))$
for any $V\in Ir\subset \Oc_0(G,q)$. Since $a(\la)$
 is rational, this means that
$a(\la)$ does not depend on $\la$.
\end{proof}

\subsection{Asymptotics of $J_{V,W}(\la)$ and $R_{V,W}(\la)$}
First assume that $q=1$ and $A=U(\g)$.
Consider $J_{V,W}(\la)$. Change variables $\la\to \la/\gm$
where $\gm\in \C^*$. Then $J_{V,W}(\la/\gm)$ has the form
$$
J_{V,W}(\la/\gm)=1 + \gm j_{V,W}(\la)+ O(\gm^2).
$$
To describe $j_{V,W}(\la)$ we fix notations. Namely, we
fix an invariant nondegenerate bilinear form  $(\cdot , \cdot)$  on
$\g$. The bilinear form identifies $\g$ and $\g^*$.
For any positive root $\al$, fix generators
$e_\al \in \g_\al, f_\al \in \g_{-\al},$ such that $h_\al=[e_\al,f_\al]$
has the property $<h_\al, \la>= (\al,\la)$ for all $\la\in \h^*$.
\begin{thm}\label{}
We have $j_{V,W}(\la)=j(\la)|_{V\T W}$, where $j(\la)\in \n_-\T \n_+$ and
\bean\label{as}
j(\la)= -\sum_{\al>0} { f_\al\T e_\al \over (\la, \al)}.
\notag
\eean
\end{thm}
\begin{corollary}
For $q=1$ and $A=U(\g)$, we have
\bean
J_{V,W}\,(u\T w)\,=\,u\T w \,-\,
\sum_{\al>0}\,{ f_\al\T e_\al \over (\la, \al)}\,u\T w + O({1\over |\la|^2})\,,
\notag
\\
R_{V,W}\,(u\T w)\,=\,u\T w \,+\,
\sum_{\al>0}\,{ f_\al\T e_\al - e_\al\T f_\al 
\over (\la, \al)}\,u\T w + O({1\over |\la|^2})\,
\notag
\eean
\end{corollary}
{\bf Proof of the Theorem.}
Let $w\in W$.
Consider the intertwining operator 
$\Phi^w_\la : M_\la\to M_{\la - wt(w)}\T W$. Let $v_\la\in M_\la$
be the highest weight vector
(we write $M_\la$ for $M_\la^+$). It follows from [ES1], Chapter 3
that 
\bean
\Phi^w_\la v_\la=v_{\la-wt(w)}\T w - \sum_{\al > 0}
{1\over (\la, \al)}\,f_\al v_{\la-wt(w)}\T e_\al w +
O({1\over |\la|^2}).
\notag
\eean
Now computing the leading term of the composition
$\Phi^u_{\la-wt(w)} \Phi^w_\la v_\la$ we conclude that
$$
J_{V,W}\,(u\T w)\,=\,u\T w \,-\,
\sum_{\al>0}\,{ f_\al\T e_\al \over (\la, \al)}\,u\T w + O({1\over |\la|^2}).
$$ 
This proves the Theorem.
$\square$

Let $q\neq 1$, $A=U_q(\g)$ and $r=$dim $\h$.
It is well known that $\Rm  \in A\hat {\T} A$ has the form
$\Rm=\Rm_0Q$ where $Q\in U_q(\h)\hat {\T}U_q(\h)$ 
is a suitable invertible element, and $\Rm_0 = 1 + $(a strictly
upper triangular element in $U_q(\n_+)\hat{\T}U_q(\n_-)$).
\begin{thm}\label{asq}
For $|q|< 1$ and $A=U_q(\g)$, 
\begin{enumerate}
\item[I.]
$J_{V,W}(\la) \to 1$, when $\la \in T=(\C^*)^r$ tends to infinity along
the positive alcove, and
$J_{V,W}(\la) \to \Rm_0^{21}$, when $\la $ tends to infinity along
the negative alcove.
\item[II.]
$R_{V,W}(\la) \to \Rm^{21}$, when $\la $ tends to infinity along
the positive alcove, and
$R_{V,W}(\la) \to Q\Rm Q^{-1}$, 
when $\la $ tends to infinity along the negative alcove.
\end{enumerate}
\end{thm}

{\bf Proof.}
It is clear that statement II follows from I. The first statement of I 
follows from \cite{ES2}, Chapter 2, as explained in the proof of Theorem 
\ref{nond2pt}. So it remains to prove the second statement of I. 

It follows from Proposition 19.3.7 in [L] that the asymptotics 
of the Shapovalov form on $M_\la=U_q(\n_-)$ 
 for $\la$ tending to $\infty$ in the negative alcove 
equals to the Drinfeld form on $U_q(\n_-)$ (i.e. the form 
which defines an injective map of Hopf algebras $U_q({\frak{b}}_+)$ to its dual). 
This fact together with the explicit formula for the intertwining 
operator via the Shapovalov form ([ES2], Chapter 2) implies 
the second statement of I. 

\subsection{Proof of part II of Theorem \ref{iso}}\label{proof}
First assume that $q=1$ and $A=U(\g)$. Let $W, U \in \Oc_0(G)$
and $b \in \Hom_{\Repq}(F(W),F(U))$. Recall that $b\in \End_\C (W,U)$
does not depend on $\la$ by Lemma \ref{const}.

\begin{lemma}
The linear operator $b$ commutes with the action of elements $e_\al, 
f_\al$ where $\al$ is any positive root.
\end{lemma}
\begin{corollary}
$b\in \Hom_{\Oc_0(G,q)}(W,U)$.
\end{corollary}
The Corollary implies part II of Theorem \ref{iso} for $q=1$.

{\bf Proof of the Lemma.}
We prove the Lemma for $W=U$. For $W\neq U$, the proof is similar.
For any $V\in Ir \subset \Oc_0(G)$, we have
$[1\T b, R_{V,W}(\la)]=0$. Hence,
$[1\T b, (R_{V,W}(\la)-1)|\la|]=0$.
Setting $\la=t\la_0$ and taking the limit $t\to \infty$, we get
$$
\sum_{\al>0}\,{ f_\al\T [b,e_\al] - e_\al\T [b,f_\al] 
\over (\la, \al)}=0.
$$
Since there exists $V$ such that the linear operators $e_\al|_V,
f_\al|_V$ are linear independent in $\End_\C(V)$, we get the Lemma.
$\square$

Part II of Theorem \ref{iso} for $|q|<1$
follows similarly from Theorem \ref{asq}.
Namely, from Theorem \ref{asq}
we get that any intertwining operator $b$ 
must commute with all elements of the form 
$(f\T 1)(\Rm^{21})$ and $(f\T 1)(Q\Rm Q^{-1})$ ($f\in U_q(\g)^*$), 
which obviously generate $U_q(\g)$. (Here if $X=\sum a_i\T b_i$ then 
$(f\T 1)(X)$ denotes $\sum f(a_i)b_i$.)
Thus, $b$ has to commute with $U_q(\g)$, 
Q.E.D. 

For $|q|>1$, the proof is analogous.

\section{Appendix: Fusion matrices and 6j-symbols}

In this appendix we discuss the relationship between 
fusion matrices introduced in Section 2, and 6j-symbols, 
for the Lie algebra $sl(2)$. For quantum $sl(2)$, the relationship 
is the same. 

Recall the definition of 6j-symbols (see e.g. [CFS], p.29). 
Let $V_a$, $a\in \Z_+/2$, 
be the irreducible representation of $sl(2)$ with spin $a$.
Let $v_a$ be the highest weight vector of $V_a$, and $v_{a,n}=f^nv_a$.  
Let $\phi_{a}^{bc}:V_a\to V_b\otimes V_c$ be the intertwiner 
such that $\phi_a^{bc}v_a=v_b\otimes v_{c,b+c-a}+l.o.t.$ (here
l.o.t. is ``lower order terms'').The 6j-symbol is defined by 
the formula 
$$
(1\otimes \phi_j^{bc})\phi_k^{aj}=\sum_n\left(\matrix a&b&n\\ c&k&j
\endmatrix\right)(\phi_n^{ab}\otimes 1)\phi_k^{nc}.
$$
The 6j-symbols not defined in this way are defined to be zero. 

{\bf Remark.} Our definition coincides with the standard one only up 
to normalization. Namely, it is more common to use a different normalization
of the operators $\phi_a^{bc}$, which results in a different 
normalization of the 6j-symbols. 

Now define $J_{bc}(\lambda):=J_{V_bV_c}(\lambda)$. 
The next proposition, which gives a connection between fusion matrices 
and 6j-symbols, follows easily from the definitions.

\begin{proposition} For any $k\in \Z_+/2$, one has 
$$
\sum_n 
\left(\matrix a&b&n\\ c&k&j
\endmatrix\right)(v_{b,b-n+a}\otimes v_{c,c-k+n})=
J_{bc}^{-1}(k)\phi_j^{bc}v_{j,j-k+a}.
$$
\end{proposition}

Thus $J_{bc}(k)$ is the unique rational function of $k$ which 
satisfies the above equation for $k\in\Z_+/2$. 

It is easy to check that under this correspondence, the 2-cocycle condition 
for $J(\lambda)$ corresponds to the Elliott-Biedenharn identity 
for 6j-symbols 
\cite{CFS} (known to mathematicians as the Maclane pentagon relation).
The dynamical Yang-Baxter equation for 
$R(\lambda)=J^{-1}(\lambda)J^{21}(\lambda)$
corresponds to the star-triangle relation. 

\section{Appendix: Recursive relations for fusion matrices}

In [A], the authors defined fusion matrices as unique solutions 
of certain linear equations, and checked that they satisfy 
the 2-cocycle condition. In this appendix, we will show that our 
fusion matrices satisfy the same linear equations, which implies 
that they are the same fusion matrices as in [A].

We will use a finite-dimensional version of the quantum Knizhnik-Zamolodchikov 
equations, which were deduced by Frenkel and Reshetikhin for quantum affine 
algebras. 
Consider the function $\Psi_{wv}(\lambda)\in \End(W\otimes V)$ given by
$$
\Psi_{wv}(\lambda):=
J_{WV}(\lambda)(w\otimes v)=<(\Phi_{\lambda-\lambda_v}^w\otimes 1)
\Phi_\lambda^v>
$$
where the notation $<,>$ was defined in Section 2.4. It follows from
a finite dimensional degeneration of 
the Frenkel-Reshetikhin theorem (Theorem 10.3.1 in [EFK]) that 
this function satisfies the following version of the 
quantum Knizhnik-Zamolodchikov 
equations:
$$
q^{2(\lambda_v,\lambda+\rho)-(\lambda_v,\lambda_v)}\Psi_{wv}(\lambda)=
R_{VW}^{21}(1\otimes q^{2\lambda-\lambda_v-\lambda_w+2\rho})\Psi_{wv}(\lambda).
$$
This implies that 
$$
J_{WV}(\lambda)(1\otimes q^{2(\lambda+\rho)-\sum x_i^2})=
R_{VW}^{21}q^{-\sum x_i\otimes  x_i}(1\otimes q^{2(\lambda+\rho)-\sum x_i^2})
J_{WV}
(\lambda)
$$
It is easy to see that the last equation is (up to simple changes 
of variable) the same as relation (18) in [A]. 

A similar computation is valid for an arbitrary 
quantized Kac-Moody algebra. This computation yields the linear relation 
for $J$ discussed in \cite{JKOS}.


\begin{thebibliography}
\normalsize

\bibitem[A] {A}
D. Arnaudon, E. Buffenoir, E. Ragoucy and Ph. Roche,
{\it Universal Solutions of Quantum Dynamical Yang-Baxter Equations},
preprint q-alg/9712037.

\bibitem[BBB] {BBB}
O.Babelon, D.Bernard, E.Billey, 
{\it A quasi-Hopf algebra interpretation of quantum 3j and 6j symbols and 
difference equations}, Phys.Lett.B, 375:89-97,1996.

\bibitem[CFS]{CFS} J.Carter, D.Flath, M.Saito,
{\it The classical and quantum 6j-symbols}, Math.Notes, 
Princeton University Press, Princeton, 1995.

\bibitem[CP]{CP} V. Chari and A. Pressley
{\it A Guide to Quantum Groups}, Cambridge University Press, 1994.




\bibitem[D]{D} V. Drinfeld,
{\it Quasi-Hopf Algebras},
Leningrad\ Math.\ J.\ 1 (1990), 1419--1457.

\bibitem[EFK]{EFK} P. Etingof, I.Frenkel, and A. Kirillov Jr, 
{\it Lectures on representation theory and Knizhnik-Zamolodchikov 
equations}, AMS, Providence, 1998. 

\bibitem[EK]{EK} P. Etingof and A. Kirillov Jr,
{\it
Representation-theoretic proof of Macdonald
inner product and symmetry identities,} Comp. Math., vol. 102, pages 179-202
1996. 


\bibitem[ES1]{ES1} P. Etingof and K. Styrkas
{\it   Algebraic integrability of Schr\"odinger
operators and representations of Lie algebras}, hep-th
9403135, Comp. Math., v. 98, No. 1, p. 91-112, 1995.



\bibitem[ES2]{ES2} P. Etingof and K. Styrkas
{\it Algebraic integrability of Macdonald operators 
and representations of quantum groups, q-alg 9603022,  to appear in
Comp. Math.} 



\bibitem[EV1]{EV1} P. Etingof and A. Varchenko
{\it Geometry and classification of solutions of the classical dynamical
Yang-Baxter equation},
Commun.\ Math.\ Phys.\  (1998), 

\bibitem[EV2]{EV2} P. Etingof and A. Varchenko
{\it Solutions of the quantum dynamical
Yang-Baxter equation and dynamical quantum groups}, preprint, 1997.
%submitted to Commun.\ Math.\ Phys.\  (1998), 





\bibitem[F]{F} G. Felder, {\it Conformal field theory and integrable
systems associated to elliptic curves},
Proceedings of the International Congress of Mathematicians,
Z\"urich 1994, p.\ 1247--1255, Birkh\"auser, 1994;
{\it Elliptic quantum groups,} preprint hep-th/9412207,
to appear in the Proceedings of the ICMP, Paris 1994.

\bibitem[FR]{FR} I. Frenkel and N. Reshetikhin,
{\it Quantum affine algebras and holonomic
difference equations},
Commun.\ Math.\ Phys.\ 146 (1992), 1--60.

\bibitem[FTV]{FTV} G. Felder, V. Tarasov and A. Varchenko,
{\it Solutions of the elliptic qKZB equations and Bethe
ansatz I,} Amer. Math. Soc. Transl. (2) Vol. 180, 1997, 45-75.

% preprint 1996, q-alg/9606005, to appear
%in the volume dedicated to V.I.Arnold's 60-th birthday.

\bibitem[FV1]{FV1} G. Felder and A. Varchenko,
{\it On representations of the elliptic quantum group
$E_{\tau,\eta}(sl_2)$}, Commun.\ Math.\ Phys.\ 181 (1996),
746--762.

\bibitem[FV2]{FV2} G. Felder and A. Varchenko,
{\it Algebraic Bethe ansatz for the elliptic quantum
group $E_{\tau,\eta}(sl_2)$}, Nuclear Physics B 480 (1996)
485-503.
%preprint q-alg/9605024

\bibitem[JKOS]{JKOS} M. Jimbo, H. Konno, S. Odake, J. Shiraishi,
{\it Quasi-Hopf twistors for elliptic quantum groups}, 
q-alg/9712029. 

\bibitem[L]{L} G.Lusztig,
{\it Introduction to quantum groups,} Birkhauser, Boston, 1993.

\bibitem[M]{M} D. Mumford, 
{\it Abelian varieties,} Tata Studies in Math., Oxford Univ. Press, 2nd ed., 
1975.


\bibitem[TV1]{TV1} V. Tarasov and A. Varchenko,
{\it Geometry of q-hypergeometric functions as a bridge between Yangians
and quantum affine algebras}, 
Invent. Math. 128, 501-588 (1997)

\bibitem[TV2]{TV2} V. Tarasov and A. Varchenko,
{\it Geometry of q-Hypergeometric Functions,
Quantum Affine Algebras and Elliptic Quantum Groups},
preprint q-alg/9703044.

\bibitem[Xu]{Xu}
{\it Quantum groupoids and universal dynamical R-matrices},
IHES preprint, June 1998. 


\end{thebibliography}
\end{document}